\documentclass[journal]{new-aiaa}
\usepackage[utf8]{inputenc}
\usepackage{textcomp}
\usepackage{graphicx}
\usepackage{amsmath}
\usepackage{subcaption}
\usepackage[version=4]{mhchem}
\usepackage{siunitx}
\usepackage{longtable,tabularx}
\usepackage{circledsteps}
\usepackage[linesnumbered,ruled,vlined]{algorithm2e}

\usepackage{booktabs}

\title{Parameter Optimization in Trajectory Planning via Differentiable Convex Programming}

\author{Ziqi Xu \footnote{Ph.D. Student, School of Astronautic; xzq24@buaa.edu.cn.}}
\affil{Beihang University, Beijing, 100191, People's Republic of China}

\author{Lin Cheng\footnote{Associate Professor, School of Astronautics; chenglin5580@buaa.edu.cn.}, Di Wu \footnote{Associate Professor, School of Astronautics, wudi2025@buaa.edu.cn, Member AIAA.} and Shengping Gong\footnote{Professor, School of Astronautics, gongsp@buaa.edu.cn, Senior Member AIAA (Corresponding Author).}}
\affil{Beihang University, Beijing, 100191, People's Republic of China}
\affil{State Key Laboratory of High-Efficiency Reusable Aerospace Transportation Technology, Beijing, 102206, People's Republic of China}

\begin{document}

\maketitle

\begin{abstract}    
    Sequential convex programming has been established as an effective framework for solving nonconvex trajectory planning problems. However, its performance is highly sensitive to problem parameters, including trajectory variables, algorithmic hyperparameters, and physical vehicle parameters. This paper introduces a differentiable sequential convex programming framework that integrates differentiable convex optimization with sequential convex programming to enable end-to-end parameter optimization. By deriving first-order sensitivity relations of second-order cone programming solutions with respect to problem data, exact gradients of trajectory performance metrics with respect to arbitrary parameters are obtained and propagated through iterations. The effectiveness of the proposed framework is validated through three representative applications: optimal terminal-time prediction for powered landing, trust-region penalty optimization in subproblems, and surface-to-mass ratio optimization for hypersonic gliding vehicles. Simulation results show that the proposed framework enables reliable gradient-based parameter learning and significantly improves numerical performance, convergence behavior, and design efficiency. These results indicate that the differentiable sequential convex programming framework provides a powerful and general tool for vehicle design, mission optimization, and hyperparameter selection in aerospace trajectory planning.
\end{abstract}

\section{Introduction}

Trajectory planning plays a central role in the performance and success of aerospace missions. Its computational efficiency and convergence reliability directly determine whether guidance and control can be executed in real time \cite{heReviewDatadrivenComputational2025,wangSurveyConvexOptimization2024}. Sequential convex programming (SCP) is an effective approach for trajectory optimization in applications such as powered landing and hypersonic vehicle entry \cite{lu2025model,liu2016entry}. Recent advances in SCP have improved convergence and runtime through enhanced constraint handling, relaxation techniques \cite{szmukSuccessiveConvexificationFuelOptimal2016}, and discretization strategies \cite{kimOptimalMidcourseGuidance2023a}. Beyond algorithmic improvements, however, a wide range of problem parameters, such as initial trajectory guesses, trust-region penalty coefficients, normalization scales, and aerodynamic shape parameters, exert substantial influence on the efficiency and success rate of SCP-based trajectory planning. The systematic and quantitative analysis of how these parameters affect convergence behavior and mission performance remains limited. Motivated by this gap, the present study integrates differentiable convex optimization with SCP to construct a differentiable sequential convex programming (DSCP) framework capable of computing end-to-end sensitivities from problem parameters to trajectory-level performance metrics, thereby enabling direct parameter optimization.


To investigate parameter influence in a structured manner, the parameters involved in SCP-based trajectory planning are categorized into three classes. 
The first class consists of the optimization variables intrinsic to the trajectory planning problem. To reduce problem dimension or enable real-time approximation, several works employ neural networks (NNs) to learn subsets of optimization variables. For example, NNs are trained to predict terminal time for powered landing from datasets generated by lossless convex programming \cite{liFreeFinalTimeFuelOptimal2022a}. In \cite{youLearningBasedOnboardGuidance2021a,YOU2022129}, key parameters of the powered landing problem are derived from the first-order necessary conditions of optimal control, and neural networks are used to approximate these parameters, which significantly reduces the computational burden of online trajectory planning. Reinforcement-learning approaches can also learn control strategies \cite{takahashiAutonomousReconnaissanceTrajectory2023, yangReinforcementLearningBasedRobustGuidance2024a}, but these methods focus on feedback decision-making rather than producing full trajectory profiles.


The second category includes algorithmic parameters used in constructing SCP subproblems, such as trust-region penalties, normalization coefficients, and constraint penalty parameters, which are significant in the convergence behavior of SCP. In particular, appropriate normalization improves the conditioning of the Karush-Kuhn-Tucker (KKT) system \cite{domahidiECOSSOCPSolver2013,kamathOptimalPreconditioningOnline2025}. Meanwhile, trust-region mechanisms are essential to ensuring the validity of linearization and the feasibility of subproblems, while also influencing the convergence rate \cite{wangSurveyConvexOptimization2024}. Existing studies have proposed various strategies, including dual-variable-based detection of trajectory oscillations for penalization \cite{baeNewTrustRegionConstraint2025}, analysis of soft trust-region SCP showing that first-order penalties may not ensure optimal convergence and the design of higher-order penalties \cite{xieHigherOrderSoftTrustRegionBasedSequential2023}, and hybrid-order trust-region formulations that balance low-order speed and high-order optimality \cite{xieHybridorderSoftTrust2024}. In \cite{fletcherNonlinearProgrammingPenalty2002}, it is pointed out that selecting appropriate soft trust-region parameters remains a major design difficulty for nonlinear programming. To address this issue, \cite{xieHighDimensionalMeritFunctionBasedSequential2025} proposes an SCP formulation based on a high-dimensional merit function. By embedding the merit function directly into the optimization process, this method removes the need for explicit selection of soft trust-region parameters and thereby alleviates the associated tuning burden. In contrast, \cite{mceowenAutotunedPrimalDual2025} introduces an alternative approach by separating nonconvex and convex constraints and developing an online trust-region adaptation scheme. This autotuned framework automatically adjusts trust-region bounds during the SCP iterations, improving convergence behavior and reducing reliance on manually chosen trust-region parameters.


The third category consists of mission-level parameters in trajectory planning or physical parameters associated with the vehicle dynamics \cite{hoburgGeometricProgrammingAircraft2012,zheng2025model}.  
In \cite{luPropellantOptimalPoweredDescent2023,luFastRobustOptimization2024}, a bilevel optimization framework is employed to determine the optimal switching point from the reentry phase to the powered-descent phase, treated as a key mission parameter. In the context of conceptual and preliminary vehicle design, aerodynamic-shape optimization typically requires extensive aerodynamic databases. Each aerodynamic model must be evaluated through comprehensive trajectory simulations and control-performance assessments. To reduce the computational cost of Computational Fluid Dynamics (CFD), existing acceleration techniques rely heavily on surrogate modeling for aerodynamic prediction. For example, \cite{jiUnsteadyAerodynamicModeling2024} used transfer learning and model fine-tuning to learn aerodynamic characteristics from CFD data, thereby accelerating flight-trajectory simulation. Furthermore, \cite{luoIntegratedOptimizationDespin2025} integrated surrogate models with multiobjective optimization to design a de-spin actuator, simultaneously addressing requirements on flight stability, range capability, and other key performance metrics.  Despite these advances, there remains a lack of methods that directly leverage trajectory planning simulations to optimize aerodynamic parameters themselves. In particular, efficient approaches that treat trajectory-performance metrics as objective functions for aerodynamic-parameter tuning are still sparsely explored.


To address the challenges associated with predicting and designing the three categories of parameters described above, the differentiable convex optimization approach adopted in this work enables the computation of exact gradients of trajectory performance metrics with respect to these parameters. This makes gradient-based optimization applicable to all three classes of parameters. Existing research in differentiable convex optimization generally derives the Jacobian of the optimal solution with respect to problem data through the KKT conditions or primal-dual residual systems.
OptNet formulated a differentiable structure for quadratic programs (QPs) and developed a parallelizable QP solver that supports efficient batched forward and backward passes \cite{amosOptNetDifferentiableOptimization2021}.
DiffCP \cite{agrawalDifferentiableConvexOptimization2019} and DiffQCP \cite{healeyDifferentiatingQuadraticCone2025} extended this idea to general cone programs by solving coupled primal-dual residual equations to obtain parameter sensitivities. These methods provide a unified treatment for differentiable conic optimization.
BPQP further exploited the KKT structure and the characterization of active sets, reformulating the backward pass as another QP 
\cite{panBPQPDifferentiableConvex}. This construction substantially improves the computational efficiency of the gradient propagation.
Based on these differentiable optimization layers, prior work has integrated differentiable model predictive control into reinforcement-learning frameworks \cite{romeroActorCriticModelPredictive2024}. This integration leverages the exploration capabilities of reinforcement learning together with the differentiability of optimization, improving policy-learning performance. Differentiable convex optimization has also been applied in  trajectory planning of unmanned aerial vehicles (UAVs) \cite{jiangSelfSupervisedLearningApproach2025}. In this line of work, a single differentiable optimization layer is embedded within a neural network to train a trajectory planning model, enabling real-time UAV trajectory generation from depth maps.

Building upon these developments, this work introduces differentiable convex optimization techniques into SCP to enable principled, gradient-based parameter optimization in trajectory planning. The key contributions are as follows:  
(1) first-order sensitivity relations between second-order cone programming (SOCP) optimal solutions and problem data are derived via dual variables and Lagrangian conditions, and embedded into SCP to construct a differentiable SCP framework that provides exact gradients of trajectory-level performance with respect to any parameter;  
(2) using the DSCP framework, three representative parameter optimization tasks are conducted: prediction of optimal terminal time for powered landing (a nonconvex variable), optimization of trust-region penalty coefficients for SCP (algorithmic parameters), and optimization of the surface-to-mass ratio of a hypersonic gliding vehicle (physical design parameter). These three cases span the major parameter categories in SCP-based trajectory planning and demonstrate the effectiveness and generality of the proposed method.

The structure of the paper is as follows. Section~\ref{ch2} summarizes the construction of SCP subproblems and clarifies the parameter classes examined in this work. Section~\ref{ch3} develops the proposed differentiable SCP framework. Section~\ref{ch4} demonstrates the parameter-optimization results for the three categories of interest. Section~\ref{ch5} provides concluding remarks. 

\section{Problem Statement}
\label{ch2}
In this section, the trajectory planning framework is presented based on SCP, including the discretization scheme, the linearization procedure, and the resulting SOCP formulation. We further clarify the categories of parameters involved in the SCP process. Consider the trajectory planning problem with time-invariant system dynamics 
\begin{equation}
    \label{equation:1}
    \begin{aligned}
    \min \quad&J=\Phi\left(\boldsymbol{x}(t_{\mathrm{f}})\right)+\int_{t_0}^{t_{\text{f}}}\Theta\left(\boldsymbol{x},\boldsymbol{u}\right)\text{d}t\\
    \text{s.t.}\quad&\dot{\boldsymbol{x}}=\boldsymbol{f}\left(\boldsymbol{x},\boldsymbol{u}\right)\\
    &\boldsymbol{h}(\boldsymbol{x},\boldsymbol{u})= \boldsymbol{0}\\
    &\boldsymbol{g}(\boldsymbol{x},\boldsymbol{u})\leq \boldsymbol{0}\\
    \end{aligned}
\end{equation}
where $\boldsymbol{x}$, $\boldsymbol{u}$ and $J$ denote the state and control variables and the performance index. To transform the continuous optimal control problem into a finite-dimensional optimization problem, all the constraints must be discretized, especially the system dynamics equation. In this work, a trapezoidal-rule-inspired discretization method is adopted, following our work in \cite{zhangReentryTrajectoryOptimization2025}. Let $N$ denote the number of discrete time intervals, each with duration $T$. The discrete form of the dynamics constraints is given by 
\begin{equation}
    \label{equation:ch2-disdy}
    \boldsymbol{x}[n+1]-\boldsymbol{x}[n]=T\boldsymbol{f}\left(\bar{\boldsymbol{x}}[n],\bar{\boldsymbol{u}}[n]\right),n\in \left[0,N - 1\right]
\end{equation}
where
\begin{equation}
    \bar{\boldsymbol{x}}[n]=\frac{\boldsymbol{x}[n+1]+\boldsymbol{x}[n]}{2},
    \bar{\boldsymbol{u}}[n]=\frac{\boldsymbol{u}[n+1]+\boldsymbol{u}[n]}{2}
\end{equation}

For both the state and control variables, the discrete dynamics constitute a nonlinear equality constraint. Within each SCP iteration, the original trajectory planning problem is linearized around a reference trajectory, thereby converting these nonlinear constraints into linear ones. Let $\boldsymbol{x}_{\text{ref}},\boldsymbol{u}_{\text{ref}},T_{\text{ref}}$ denote the reference state, control, and time, respectively. By applying a first-order Taylor expansion to the discrete dynamics in Eq.\eqref{equation:ch2-disdy}, the linearized form is obtained as follows:
\begin{equation}
    \label{equation:ch2-4}
     \delta \boldsymbol{x}[n+1] - \delta \boldsymbol{x}[n] + \boldsymbol{x}_{\text{ref}}[n+1]- \boldsymbol{x}_{\text{ref}}[n] = T_{\text{ref}} \boldsymbol{f}_{\text{ref}}+T_{\text{ref}}\left.\frac{\partial \boldsymbol{f}}{\partial \boldsymbol{x}}\right|_{\text{ref}}\delta \bar {\boldsymbol{x}}[n] +T_{\text{ref}}\left.\frac{\partial \boldsymbol{f}}{\partial \boldsymbol{u}}\right|_{\text{ref}}\delta \bar {\boldsymbol{u}}[n] + \delta T \boldsymbol{f}_{\text{ref}}
\end{equation}
where 
\begin{equation}
    \label{equation:ch2-5}
    \begin{aligned}
      & \left.\frac{\partial \boldsymbol{f}}{\partial \boldsymbol{x}}\right|_{\text{ref}}=\left.\frac{\partial \boldsymbol{f}}{\partial \boldsymbol{x}}\right|_{
        \bar {\boldsymbol{x}}_{\text{ref}},\bar {\boldsymbol{u}}_{\text{ref}}},\left.\frac{\partial \boldsymbol{f}}{\partial \boldsymbol{u}}\right|_{\text{ref}}=\left.\frac{\partial \boldsymbol{f}}{\partial \boldsymbol{u}}\right|_{
        \bar {\boldsymbol{x}}_{\text{ref}},\bar {\boldsymbol{u}}_{\text{ref}}},
        \boldsymbol{f}_{\text{ref}} = \boldsymbol{f}(\bar{\boldsymbol{x}}_{\text{ref}},\bar{\boldsymbol{u}}_{\text{ref}}) \\ &
    \delta\bar{\boldsymbol{x}}[n]=\frac{\delta\boldsymbol{x}[n+1]+\delta\boldsymbol{x}[n]}{2},
    \delta\bar{\boldsymbol{u}}[n]=\frac{\delta\boldsymbol{u}[n+1]+\delta\boldsymbol{u}[n]}{2},n\in \left[0,N - 1\right] \\
    & \delta\boldsymbol{x}=\boldsymbol{x}-\boldsymbol{x}_{\mathrm{ref}},\delta\boldsymbol{u}=\boldsymbol{u}-\boldsymbol{u}_{\mathrm{ref}},\delta T=T-T_{\mathrm{ref}}
    \end{aligned}
\end{equation}

By combining Eq.\eqref{equation:ch2-4} and Eq.\eqref{equation:ch2-5}, it yields
\begin{equation}
    \label{equation:6}
    \begin{aligned}
        \boldsymbol{d}_{\text{ref}} &= \left(\frac{T_{\text{ref}}}{2}\left.\frac{\partial \boldsymbol{f}}{\partial \boldsymbol{x}}\right|_{\text{ref}} + \boldsymbol{I}\right)\delta \boldsymbol{x}[n] + \left(\frac{T_{\text{ref}}}{2}\left.\frac{\partial \boldsymbol{f}}{\partial \boldsymbol{x}}\right|_{\text{ref}} - \boldsymbol{I}\right)\delta \boldsymbol{x}[n+1] \\ &+ \frac{T_{\text{ref}}}{2}\left.\frac{\partial \boldsymbol{f}}{\partial \boldsymbol{u}}\right|_{\text{ref}}\delta  \boldsymbol{u}[n] + \frac{T_{\text{ref}}}{2}\left.\frac{\partial \boldsymbol{f}}{\partial \boldsymbol{u}}\right|_{\text{ref}} \delta \boldsymbol{u}[n+1]+\boldsymbol{f}_{\text{ref}} \delta T
    \end{aligned}
\end{equation}
where $\boldsymbol{d}_{\text{ref}}=(\boldsymbol{x}_{\text{ref}}[n+1]- \boldsymbol{x}_{\text{ref}}[n])-T_{\text{ref}} \boldsymbol{f}_{\text{ref}}$ denotes the discretization residual associated with the reference trajectory. This term mitigates the adverse effects that may arise in subsequent iterations when the reference trajectory does not fully satisfy the system dynamics. Using the same procedure as for the dynamics constraints, the cost function and path constraints in Eq.\eqref{equation:1} are discretized and linearized into the following form:
\begin{itemize}
    \item inequality constraints
    \begin{equation}
    \label{equation:7}
        \boldsymbol{g}\left(\boldsymbol{x}_{\text{ref}}[n],\boldsymbol{u}_{\text{ref}}[n]\right)+\left.\frac{\partial \boldsymbol{g}}{\partial \boldsymbol{x}}\right|_{\boldsymbol{x}_{\text{ref}}[n],\boldsymbol{u}_{\text{ref}}[n]}\delta \boldsymbol{x}[n] +\left.\frac{\partial \boldsymbol{g}}{\partial \boldsymbol{u}}\right|_{\boldsymbol{x}_{\text{ref}}[n],\boldsymbol{u}_{\text{ref}}[n]}\delta \boldsymbol{u}[n] \leq 0
        \end{equation}
    \item equality constraints 
    \begin{equation}
    \label{equation:8}
        \boldsymbol{h}\left(\boldsymbol{x}_{\text{ref}}[n],\boldsymbol{u}_{\text{ref}}[n]\right)+\left.\frac{\partial \boldsymbol{h}}{\partial \boldsymbol{x}}\right|_{\boldsymbol{x}_{\text{ref}}[n],\boldsymbol{u}_{\text{ref}}[n]}\delta \boldsymbol{x}[n] +\left.\frac{\partial \boldsymbol{h}}{\partial \boldsymbol{u}}\right|_{\boldsymbol{x}_{\text{ref}}[n],\boldsymbol{u}_{\text{ref}}[n]}\delta \boldsymbol{u}[n] = 0
    \end{equation}
    \item the cost function 
    \begin{equation}
        \begin{aligned}
    J_1&=J\left(\boldsymbol{x}_{\text{ref}},\boldsymbol{u}_{\text{ref}}\right)+\left.\frac{\partial \Phi}{\partial \boldsymbol{x}}\right|_{\boldsymbol{x}_\text{ref}[N]}\delta \boldsymbol{x}[N]+T\sum_{n=0}^{N}
    \delta \Theta[n] \text{d}t +\delta T\sum_{n=0}^{N}\Theta[n]\\
        \delta \Theta&=\left.\frac{\partial \Theta}{\partial \boldsymbol{x}}\right|_{\boldsymbol{x}_{\text{ref}},\boldsymbol{u}_{\text{ref}}}\delta \boldsymbol{x} +\left.\frac{\partial \Theta}{\partial \boldsymbol{u}}\right|_{\boldsymbol{x}_{\text{ref}},\boldsymbol{u}_{\text{ref}}}\delta \boldsymbol{u}
        \end{aligned}
    \end{equation}
\end{itemize}

The linearized cost and constraint expressions derived above yield a linear programming subproblem. However, within the SCP framework, the discretized and linearized subproblem always incorporates either trust-region constraints or trust-region penalty terms. These mechanisms restrict the allowable variations in the state and control trajectories, thereby ensuring that each update remains within a neighborhood where the linearization of the dynamics is reliable. In this work, a trust-region penalty is adopted, expressed as 
\begin{equation}
    J_{\text{trust}} = \sum_{n=0}^{N}\left\|D(\boldsymbol{\omega}_{\text{trust},\boldsymbol{x}})\delta\boldsymbol{x}[n]+D(\boldsymbol{\omega}_{\text{trust},\boldsymbol{u}})\delta\boldsymbol{u}[n]\right\|_2^2 
\end{equation}
where $D(\cdot)$ denotes a diagonal matrix formed from the vector $(\cdot)$. Unlike conventional trust-region penalty formulations, the proposed scheme assigns distinct penalty weights to different components of the state and control vectors. Although normalization is applied to nondimensionalize the optimization variables, it is generally not possible to scale all dimensions to comparable magnitudes. Adjusting the trust-region coefficients for each dimension therefore plays a critical role in influencing the convergence behavior of the SCP procedure, as demonstrated in Sec.\ref{ch4-2}. With the trust-region term incorporated, the SCP subproblem at each iteration takes the following form:
\begin{equation}
    \label{equation:P1}
    \begin{aligned}
    \min &\quad J_2=J_1+J_{\mathrm{trust}}\\
\text{s.t.}&\quad\mathrm{Eq.}\eqref{equation:6},\mathrm{Eq.}\eqref{equation:7},\mathrm{Eq.}\eqref{equation:8}\\
    \end{aligned}
\end{equation}

The above formulation constitutes a quadratic programming subproblem. Let $\boldsymbol{z}$ and $\boldsymbol{z}_{\text{ref}}$ denote the variables of the SCP subproblem and the corresponding reference trajectory. For problems with fixed terminal time, the vector typically takes the form $\boldsymbol{z}=\left\{(\boldsymbol{x}[n],\boldsymbol{u}[n]),n\in [0, N]\right\}$, whereas for problems with free terminal time, it is given by $\boldsymbol{z}=\left\{T,(\boldsymbol{x}[n],\boldsymbol{u}[n]),n\in [0, N]\right\}$. It is significant to note that, in the above SCP implementations, the variables at each iteration are the increments $\delta \boldsymbol{z}$ relative to the current reference trajectory. The updated trajectory is then obtained via $\boldsymbol{z}_{\star}=\delta\boldsymbol{z}+\boldsymbol{z}_{\mathrm{ref}}$. For notational simplicity, these two steps are combined, and $\boldsymbol{z}$ is used directly as the optimization variable throughout the presentation. 

In many trajectory planning applications, certain constraints such as thrust constraints in planetary landing are intrinsically second-order cone (SOC) constraints that should not be linearized. To accommodate this structure in a general manner, SOC constraints are directly incorporated into Eq.\eqref{equation:P1}, and then it is rewritten in the following SOCP form:
\begin{equation}
    \label{equation:P2}
    \begin{aligned}
        \min_{\boldsymbol{z}}\quad&\frac{1}{2}\boldsymbol{z}^{\top}\boldsymbol{Q}(\boldsymbol{\theta})\boldsymbol{z}+\boldsymbol{c}^{\top}(\boldsymbol{\theta})\boldsymbol{z}\\
        \mathrm{s.t.~}\quad&\boldsymbol{A}(\boldsymbol{\theta})\boldsymbol{z}=\boldsymbol{b}_0(\boldsymbol{\theta})\\
        &\boldsymbol{G}_0(\boldsymbol{\theta})\boldsymbol{z}\leq \boldsymbol{h}_0(\boldsymbol{\theta})\\
        &\boldsymbol{g_i} \coloneq \left[ g_{i,0}\left(\boldsymbol{z},\boldsymbol{\theta}\right),\boldsymbol{g}_{i,1}^\top\left(\boldsymbol{z}, \boldsymbol{\theta}\right)\right]^\top\in \mathcal{Q}_{m_i},i=1,...,N_{\mathrm{SOC}}
    \end{aligned}
    \end{equation}
where 
\begin{equation}
    \begin{aligned}
        g_{i, 0}&=\boldsymbol{a}_i^\top(\boldsymbol{\theta})\boldsymbol{z}+b_i(\boldsymbol{\theta}) \\
\boldsymbol{g}_{i,1}&=\boldsymbol{G}_i(\boldsymbol{\theta})\boldsymbol{z}+\boldsymbol{h}_i(\boldsymbol{\theta}) \\
        \mathcal{Q}_{m_i}&=\left\{(g_{i,0},\boldsymbol{g}_{i,1})|g_{i,0}\geq\left\|\boldsymbol{g}_{i,1}\right\|_2\right\}\\
    \end{aligned}
\end{equation}

$N_{\mathrm{SOC}}$ denotes the number of second-order cone constraints, and 
$m_i$ represents the dimension of the $i$-th cone. The vector $\boldsymbol{\theta}$ aggregates all parameters involved in the SCP subproblem, including the reference trajectory $\boldsymbol{z}_{\text{ref}}$, the SCP-specific parameters $\boldsymbol{\theta}_{\text{SCP}}$ and the trajectory planning parameters $\boldsymbol{\theta}_{\text{TP}}$. By repeatedly solving the subproblem in Eq.\eqref{equation:P2} and updating the reference trajectory, SCP ultimately converges to the optimal solution of the problem in Eq.\eqref{equation:1}. Although SCP is widely recognized as an effective method for trajectory planning, its convergence efficiency is highly sensitive to parameter selection, which has attracted increasing research attention in recent years \cite{xieHigherOrderSoftTrustRegionBasedSequential2023,xieHybridorderSoftTrust2024}. In addition, fixing part of the nonconvex decision variables prior to optimization, thereby treating them as parameters during the solution process, has proven to be an effective means of accelerating computation \cite{liFreeFinalTimeFuelOptimal2022a,doi:10.2514/6.2024-1760}. Motivated by these observations, this work investigates the influence of parameter optimization on the SCP-based trajectory planning. A differentiable convex optimization framework is integrated into the SCP scheme, enabling end-to-end learning of parameters through a differentiable SCP pipeline. 

\section{Differentiable Sequential Convex Programming}
\label{ch3}
This section first derives the total differential relationship between the optimal value of a single convex optimization problem and its parameters. It then presents the forward computation and backward gradient propagation procedures of the differentiable SCP (DSCP) framework, enabling end-to-end optimization of SCP parameters.
\subsection{Differentiable Second-Order Cone Programming}
The proposed method incorporates differentiable SOCP into the SCP framework and introduces several modifications tailored to the characteristics of SCP-based trajectory planning. Let $\boldsymbol{\nu}$, $\boldsymbol{\mu}$ and $\boldsymbol{\lambda}_{i},i\in [1,N_{\mathrm{SOC}}]$ denote the dual variables associated with the linear equality constraints, linear inequality constraints, and SOC constraints, respectively. The Lagrangian of the optimization problem in Eq.\eqref{equation:P2} is defined as
\begin{equation}
    L(\boldsymbol{z},\boldsymbol{\nu},\boldsymbol{\mu},\boldsymbol{\lambda})=\frac{1}{2}\boldsymbol{z}^\top \boldsymbol{Q}\boldsymbol{z}+\boldsymbol{c}^\top \boldsymbol{z}+\boldsymbol{\nu}^\top(\boldsymbol{A}\boldsymbol{z}-\boldsymbol{b})+\boldsymbol{\mu}^\top\left(\boldsymbol{G}_0\boldsymbol{z}-\boldsymbol{h}_0\right)-\sum_{i=1}^{N_{\mathrm{SOC}}}\boldsymbol{\lambda}_i^\top\boldsymbol{g}_i
\end{equation}
where $\boldsymbol{\lambda}_i=[\lambda_{i,0},\boldsymbol{\lambda}_{i,1}^\top]^\top,\lambda_{i,0}\in\mathbb{R},\boldsymbol{\lambda}_{i,1}\in \mathbb{R}^{m_i-1} $. Define $z_{\star}$ as the optimal primal solution of SOCP. The stationarity condition at the optimum is 
\begin{equation}
        \left.\frac{\partial L}{\partial \boldsymbol{z}} \right|_{\boldsymbol{z}=\boldsymbol{z}_\star}= \boldsymbol{Q}\boldsymbol{z}_{\star}+\boldsymbol{c}+\boldsymbol{A}^\top \boldsymbol{\nu}+\boldsymbol{G}_0^\top \boldsymbol{\mu}-\sum_{i=1}^{K}\left(\lambda_{i,0}\boldsymbol{a}_i+\boldsymbol{G}_{i}^{\top}\boldsymbol{\lambda}_{i,1} \right)=0
\end{equation}

The primal feasibility conditions are 
\begin{equation}
    \boldsymbol{A}\boldsymbol{z}_{\star}-\boldsymbol{b}=0,\boldsymbol{G}_0\boldsymbol{z}_{\star}-\boldsymbol{h}_0\leq0,\boldsymbol{g_i}\in \mathcal{Q}_{m_i}
\end{equation}

The dual feasibility conditions are 
\begin{equation}
    \boldsymbol{\mu} \geq 0, \boldsymbol{\lambda}_i \in \mathcal{Q}_{m_i}
\end{equation}

The complementary slackness conditions are 
\begin{equation}
    \boldsymbol{\mu}\odot  (\boldsymbol{G}_0\boldsymbol{z}_{\star}-\boldsymbol{h}_0) = 0
\end{equation}
\begin{equation}
    \langle \boldsymbol{\lambda}_i, \boldsymbol{g}_i\rangle = 
    \lambda_{i,0}g_{i,0}+
    \boldsymbol{\lambda}_{i,1}^\top \boldsymbol{g}_{i,1}= 0 
\end{equation} 

From dual feasibility and complementary slackness, it follows that if an inequality constraint is strictly satisfied at the optimum, its corresponding dual variables must be zero. Conversely, a nonzero dual variable can only arise from a constraint that is active at the optimal solution. Therefore, define the active sets for the linear inequality constraints and the SOC constraints as 
\begin{equation}
    \label{equation:active}
    \begin{aligned}
        &\mathcal A_{\mathrm{lin}} \coloneq \{ j \mid (G_0 z_{\star} - h_0)_j = 0 \}= \{j \mid \mu_{j} > 0 \,\}\\
 &\mathcal A_{\mathrm{soc}}\coloneq \{ i \mid g_{i,0}=\left\|\boldsymbol{g}_{i,1}\right\|_2 \}= \{i \mid \lambda_{i, 0} > 0 \}.
    \end{aligned}
\end{equation}
Constraints with indices $j\notin\mathcal A_{\mathrm{lin}}$ and  $i\notin\mathcal A_{\mathrm{soc}}$ do not contribute to the stationarity condition or subsequent derivations \cite{panBPQPDifferentiableConvex}. In the remainder of this section, all summations are taken only over active constraints. Based on complementary slackness, one obtains 
\begin{equation}
    \label{equation:20}
    \lambda_{i,0}g_{i,0}=-\boldsymbol{\lambda}_{i,1}^\top \boldsymbol{g}_{i,1} \leq \left|\boldsymbol{\lambda}_{i,1}^\top \boldsymbol{g}_{i,1}\right| \leq \left\|\boldsymbol{\lambda}_{i,1}\right\|_2\left\|\boldsymbol{g}_{i,1}\right\|_2,i\in \mathcal A_{\mathrm{lin}}
\end{equation}
where the last inequality follows from the Cauchy-Schwarz inequality. In addition, since both $\boldsymbol{\lambda}_i$ and $\boldsymbol{g}_i$ belong to the cone $\mathcal{Q}_{m_i}$, it holds that
\begin{equation}
    \label{equation:21}
    \lambda_{i,0} \geq \left\|\boldsymbol{\lambda}_{i,1}\right\|_2,
    g_{i,0} \geq \left\|\boldsymbol{g}_{i,1}\right\|_2 \rightarrow 
    \lambda_{i,0} g_{i,0} \geq \left\|\boldsymbol{\lambda}_{i,1}\right\|_2\left\|\boldsymbol{g}_{i,1}\right\|_2
\end{equation}
Combining Eq.\eqref{equation:20} and Eq.\eqref{equation:21}, the primal and dual cone variables satisfy the property 
\begin{equation}
    \lambda_{i,0}=\left\|\boldsymbol{\lambda}_{i,1}\right\|_2,
    g_{i,0}=\left\|\boldsymbol{g}_{i,1}\right\|_2,-\boldsymbol{\lambda}_{i,1}^\top \boldsymbol{g}_{i,1}=\left\|\boldsymbol{\lambda}_{i,1}\right\|_2\left\|\boldsymbol{g}_{i,1}\right\|_2
\end{equation}
which implies that $\boldsymbol{\lambda}_{i,1}=-\lambda_{i,0}\boldsymbol{g}_{i,1}/g_{i,0}$. Substituting this expression into the stationarity condition and applying the definition of the active set yield 
\begin{equation}
    \boldsymbol{Q}\boldsymbol{z}_{\star}+\boldsymbol{c}+\boldsymbol{A}^\top \boldsymbol{\nu}+\sum_{j\in \mathcal{A}_{\mathrm{lin}}}\mu_{j}(\boldsymbol{G}_{0})_{j}^\top+\sum_{i\in \mathcal{A}_{\mathrm{soc}}}\lambda_{i,0}\left(\boldsymbol{G}_{i}^\top\frac{\boldsymbol{g}_{i,1}}{g_{i,0}}-\boldsymbol{a}_i\right)=0
\end{equation}


Taking the total differential of the stationarity condition gives
\begin{equation}
    \begin{aligned}
&\text{d}\boldsymbol{Q}\boldsymbol{z}_{\star}+\boldsymbol{Q}\text{d}\boldsymbol{z}_{\star}+\text{d}\boldsymbol{c}+\text{d}\boldsymbol{A}^\top\boldsymbol{\nu}+\boldsymbol{A}^\top\text{d}\boldsymbol{\nu}\\+&\sum_{j\in \mathcal{A}_{\mathrm{lin}}}\left(\mathrm{d}\mu_j\left(\boldsymbol{G}_0\right)_j^\top+\mu_j\left(\mathrm{d}\boldsymbol{G}_0\right)_j^\top\right)+\sum_{i\in \mathcal{A}_{\mathrm{soc}}}\left(\mathrm{d}\lambda_{i,0}\boldsymbol{s}_{i}+\lambda_{i,0}\mathrm{d}\boldsymbol{s}_{i}\right)= 0
    \end{aligned}
\end{equation}
where
\begin{equation}
    \label{equation:26}
    \begin{aligned}
&\boldsymbol{s}_i=\boldsymbol{G}_{i}^\top\frac{\boldsymbol{g}_{i,1}}{g_{i,0}}-\boldsymbol{a}_i\\
&\mathrm{d}\boldsymbol{s}_i=(\mathrm{d}\boldsymbol{G}_i^\top)\frac{\boldsymbol{g}_{i,1}}{g_{i,0}}+\frac{\boldsymbol{G}_i^\top}{g_{i,0}}\mathrm{d}\boldsymbol{g}_{i,1}-\frac{\boldsymbol{G}_i^\top\boldsymbol{g}_{i,1}}{g_{i,0}^2}\mathrm{d}g_{i,0}-\mathrm{d}\boldsymbol{a}_i\\
&\mathrm{d}g_{i,0}=(\mathrm{d}\boldsymbol{a}_i)^\top\boldsymbol{z}_{\star}+\boldsymbol{a}_i^\top\mathrm{d}\boldsymbol{z}_{\star}+\mathrm{d}b_i\\
&\mathrm{d}\boldsymbol{g}_{i,1}=(\mathrm{d}\boldsymbol{G}_i)\boldsymbol{z}_{\star}+\boldsymbol{G}_i\mathrm{d}\boldsymbol{z}_{\star}+\mathrm{d}\boldsymbol{h}_i
    \end{aligned}
\end{equation}
and $\boldsymbol{s}_i$ denotes the outward normal direction of the $i$-th SOC at $z_{\star}$. Rearranging terms yields 
\begin{equation}
    \label{equation:27}
    \begin{aligned}
&\boldsymbol{R}\mathrm{d}\boldsymbol{z}_{\star}+\boldsymbol{A}^\top\mathrm{d}\boldsymbol{\nu}+\sum_{j\in \mathcal{A}_{\mathrm{cal}}}\mathrm{d}\mu_j\left(\boldsymbol{G}_0\right)_j^{\top}+\sum_{i\in \mathcal{A}_\mathrm{soc}}\mathrm{d}\lambda_{i,0}\boldsymbol{s}_{i}=-\boldsymbol{r}_{\boldsymbol{z}}\\
\boldsymbol{R}&=\boldsymbol{Q}+\sum_{i\in\mathcal{A}_{\mathrm{soc}}}\lambda_{i,0}\left(\frac{\boldsymbol{G}_i^\top\boldsymbol{G}_i}{g_{i,0}}-\frac{\boldsymbol{G}_i^\top\boldsymbol{g}_{i,1}\boldsymbol{a}_i^\top}{g_{i,0}^2}\right)\\
\boldsymbol{r}_{\boldsymbol{z}}&=(\mathrm{d}\boldsymbol{Q})\boldsymbol{z}_{\star}+\mathrm{d}\boldsymbol{c}+(\mathrm{d}\boldsymbol{A}^{\top})\boldsymbol{\nu}+\sum_{j\in\mathcal{A}_{\mathrm{lin}}}\mu_{j} (\mathrm{d}\boldsymbol{G}_{0})_{j}^{\top}\\
&+\sum_{i\in\mathcal{A}_{\mathrm{soc}}}\lambda_{i,0}\left[(\mathrm{d}\boldsymbol{G}_{i}^{\top})\frac{\boldsymbol{g}_{i,1}}{g_{i,0}}+\frac{\boldsymbol{G}_{i}^{\top}}{g_{i,0}}\big((\mathrm{d}\boldsymbol{G}_{i})\boldsymbol{z}_{\star}+\mathrm{d}\boldsymbol{h}_{i}\big)\right.\\
&- \frac{\boldsymbol{G}_{i}^{\top}\boldsymbol{g}_{i,1}}{g_{i,0}^{2}}\big((\mathrm{d}\boldsymbol{a}_{i})^{\top}\boldsymbol{z}_{\star}+\mathrm{d}b_{i}\big)-\mathrm{d}\boldsymbol{a}_{i}\bigg]
    \end{aligned}
\end{equation}
For the linear equality constraints and the active inequality constraints, the total differentials are 
\begin{equation}
        \boldsymbol{A}\mathrm{d}\boldsymbol{z}_{\star}=\mathrm{d}\boldsymbol{b} - \mathrm{d}\boldsymbol{A}\boldsymbol{z}_{\star}=\boldsymbol{r}_{\boldsymbol{A}}
\end{equation}
\begin{equation}
    \label{equation:G1}
\left(\boldsymbol{G}_0\right)_j\mathrm{d}\boldsymbol{z}_{\star}=\mathrm{d}\boldsymbol{h}_j-\left(\mathrm{d}\boldsymbol{G}_0\right)_j\boldsymbol{z}_{\star},j\in\mathcal{A}_{\mathrm{lin}}
\end{equation}
\begin{equation}
    \label{equation:S1}
\boldsymbol{s}_{i}^{\top}\mathrm{d}\boldsymbol{z}_{\star}=\mathrm{d}\boldsymbol{b}_i+\mathrm{d}\boldsymbol{a}_i^{\top}\boldsymbol{z}_{\star}-\frac{\boldsymbol{g}_{i,1}^{\top}}{g_{i,0}}\left(\mathrm{d}\boldsymbol{h}_i+\mathrm{d}\boldsymbol{G}_i\boldsymbol{z}_{\star}\right),i\in\mathcal{A}_{\mathrm{soc}}
\end{equation}
Concatenating Eq.\eqref{equation:G1} and Eq.\eqref{equation:S1} columnwise yields the matrix expressions that collect all active inequality constraints.  
\begin{equation}
    \label{equation:GS}
    \begin{aligned}
        &\boldsymbol{G}_{0+}\mathrm{d}\boldsymbol{z}_{\star}=\boldsymbol{r}_{\boldsymbol{G}}\\
        &\boldsymbol{S}_{+}\mathrm{d}\boldsymbol{z}_{\star}=\boldsymbol{r}_{\boldsymbol{S}}\\
    \end{aligned}
\end{equation}
with corresponding dual differentials $\mathrm{d}\boldsymbol{\mu}_{+}$ and $\mathrm{d}\boldsymbol{\lambda}_{0+}$. 

Substituting Eq.\eqref{equation:S1} into Eq.\eqref{equation:27} and cancelling common terms yields
\begin{equation}
    \label{equation:32}
    \begin{aligned}
\boldsymbol{R}\leftarrow &\, \boldsymbol{Q}+\sum_{i\in\mathcal{A}_{\mathrm{soc}}}\lambda_{i,0}\left(\frac{\boldsymbol{G}_i^\top\boldsymbol{G}_i}{g_{i,0}}-\frac{\boldsymbol{G}_i^\top\boldsymbol{g}_{i,1}\boldsymbol{g}_{i,1}^\top\boldsymbol{G}_i}{g_{i,0}^3}\right)\\
\boldsymbol{r}_{\boldsymbol{z}}\leftarrow&\, (\mathrm{d}\boldsymbol{Q})\boldsymbol{z}_{\star}+\mathrm{d}\boldsymbol{c}+(\mathrm{d}\boldsymbol{A}^{\top})\boldsymbol{\nu}+\sum_{j\in\mathcal{A}_{\mathrm{lin}}}\mu_{j} (\mathrm{d}\boldsymbol{G}_{0})_{j}^{\top}\\
&+\sum_{i\in\mathcal{A}_{\mathrm{soc}}}\lambda_{i,0}\left[(\mathrm{d}\boldsymbol{G}_{i}^{\top})\frac{\boldsymbol{g}_{i,1}}{g_{i,0}}+\frac{\boldsymbol{G}_{i}^{\top}}{g_{i,0}}\left(\boldsymbol{I}-\frac{\boldsymbol{g}_{i,1}\boldsymbol{g}_{i,1}^{\top}}{g_{i,0}^2}\right)\big((\mathrm{d}\boldsymbol{G}_{i})\boldsymbol{z}_{\star}+\mathrm{d}\boldsymbol{h}_{i}\big)\right.-\mathrm{d}\boldsymbol{a}_{i}\bigg]
    \end{aligned}
\end{equation}

Combining Eq.\eqref{equation:27} through Eq.\eqref{equation:32}, it yields the linear system  
\begin{equation}
\label{equation:KKT1}    \begin{bmatrix}\boldsymbol{R}&\boldsymbol{A}^\top&\boldsymbol{G}_{0+}^\top&\boldsymbol{S}_{+}\\\boldsymbol{A}&0&0&0\\\boldsymbol{G}_{0+}&0&0&0\\\boldsymbol{S}_{+}^\top&0&0&0\end{bmatrix}\begin{bmatrix}\mathrm{d}\boldsymbol{z}_{\star}\\\text{d}\boldsymbol{\nu}\\\text{d}\boldsymbol{\mu}_{+}\\\text{d}\boldsymbol{\lambda}_{0+}\end{bmatrix}=\begin{bmatrix}-\boldsymbol{r}_{\boldsymbol{z}}\\\boldsymbol{r}_{\boldsymbol{A}}\\\boldsymbol{r}_{\boldsymbol{G}}\\\boldsymbol{r}_{\boldsymbol{S}}\\\end{bmatrix}
\end{equation}
which is written compactly as
\begin{equation}
    \boldsymbol{H}\boldsymbol{\Delta}=\boldsymbol{r}_{\Delta}
\end{equation}
Since $\boldsymbol{R}$ is a symmetric matrix, and then the matrix $\boldsymbol{H}$ is  symmetric as well. All components of $\boldsymbol{R},\boldsymbol{A},\boldsymbol{G}_{0+},\boldsymbol{S}_{+}$ are computed at the optimal point from $\boldsymbol{z}_{\star}$ and the problem parameters. Therefore,  Eq.\eqref{equation:KKT1} implicitly characterizes the sensitivity of the optimal solution with respect to all problem parameters. Since matrices and vectors, such as $\boldsymbol{Q}(\boldsymbol{\theta})$ and $\boldsymbol{A}(\boldsymbol{\theta})$, are analytic functions of the parameter set $\boldsymbol{\theta}$, their differentials $\mathrm{d}\boldsymbol{Q},\mathrm{d}\boldsymbol{A}$ and $\mathrm{d}\boldsymbol{b}$ respect to $\boldsymbol{\theta}$ can be computed directly. Substituting these differentials into the right-hand side of Eq.\eqref{equation:KKT1} yields the partial derivatives or directional derivatives of the optimal solution with respect to the parameters $\boldsymbol{\theta}$. 

Serving as a classical sensitivity analysis tool, the formulation characterizes the first-order response of the optimal solution to perturbations in constraint boundaries, cost weights, and system parameters, thereby offering insight into how optimal trajectories and control actions vary under changes in problem data. This establishes the local linear sensitivity of parameterized SOCP solutions. At the same time, Eq.\eqref{equation:KKT1} provides an exact gradient propagation mechanism for differentiable convex optimization, enabling SOCP subproblems to be embedded seamlessly within differentiable computational pipelines. 

\subsection{Differentiable Sequential Convex Programming}

\begin{figure*}[t]
    \centering
    \includegraphics[width=\textwidth]{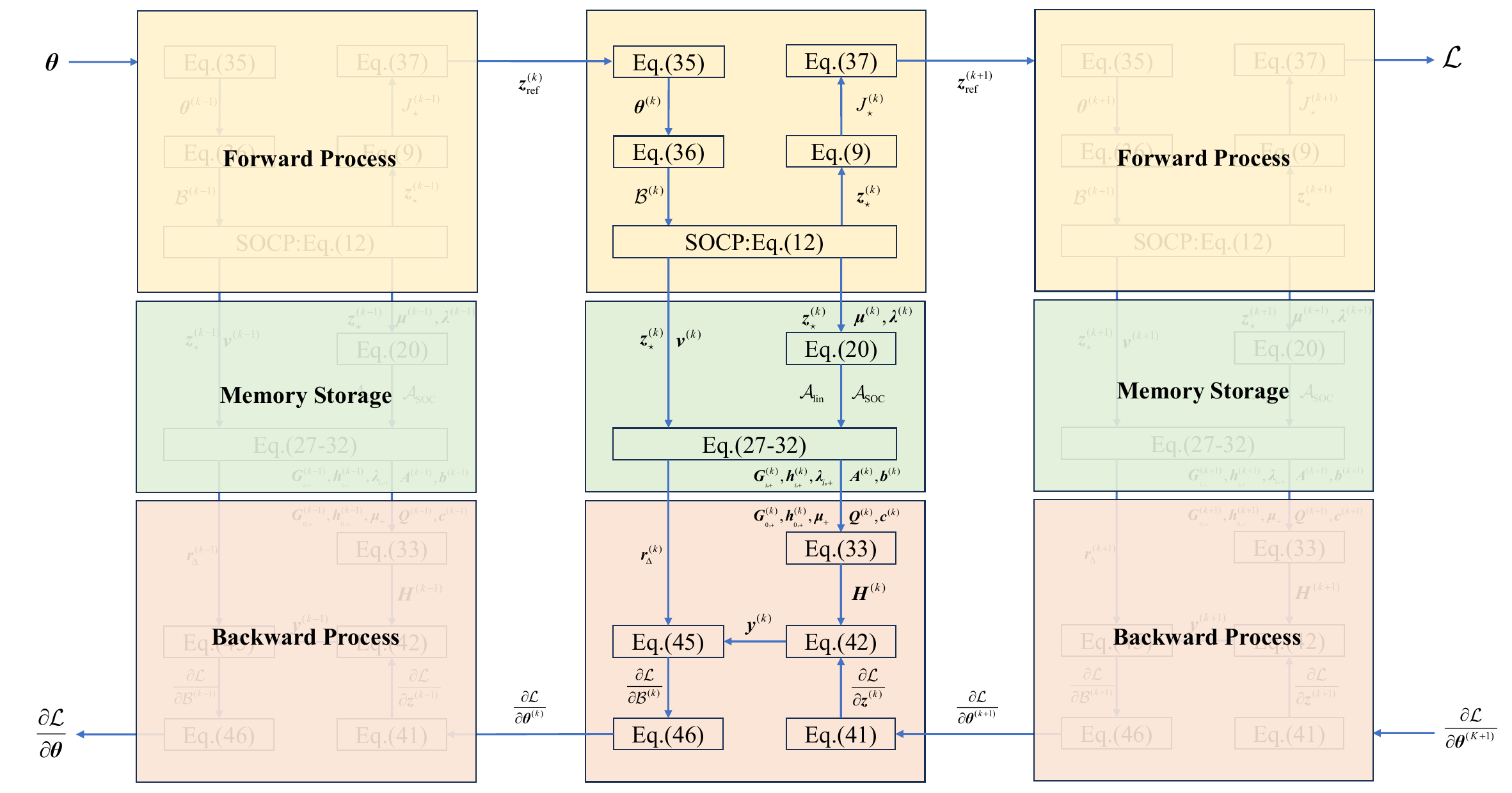}
    \caption{Differentiable Sequential Convex Programming Architecture.}
    \label{figure:algo1}
\end{figure*}

The previous subsection established the total differential relationship between the optimal solution of an SOCP and its problem parameters. Treating this mapping as an optimization layer with well-defined forward and backward computational flows, multiple such layers can be stacked to construct a fully differentiable SCP pipeline. The resulting DSCP architecture is illustrated in Fig. \ref{figure:algo1}. The overall procedure consists of three computational stages: a forward pass that executes the SCP iterations, a backward pass that propagates gradients through the entire iterative process, and a sensitivity-construction stage in which the active constraint sets and the associated sensitivity linear systems, derived from the original problem data, are assembled and stored in memory for the backward computation.
\subsubsection{Forward Process of Differentiable Sequential Convex Programming}

The forward process largely mirrors the standard SCP solution procedure. For brevity, only the elements relevant to the subsequent backward gradient computation are highlighted. Let $k$ denote the iteration index of the forward SCP process and $\boldsymbol{z}_{\mathrm{ref}}^{(k)}$ represent the reference trajectory at iteration $k$. Together with $\boldsymbol{\theta}_{\mathrm{SCP}},\boldsymbol{\theta}_{\mathrm{TP}}$, these quantities form the parameter set $\boldsymbol{\theta}^{(k)}$ as 
\begin{equation}
    \boldsymbol{\theta}^{(k)}=[\boldsymbol{z}_{\mathrm{ref}}^{(k)\top},\boldsymbol{\theta}_{\mathrm{SCP}}^{\top},\boldsymbol{\theta}_{\mathrm{TP}}^{\top}]^\top
\end{equation}
which determines the set $\mathcal{B}^{(k)}$ of all coefficient matrices and vectors of the SOCP constructed at iteration $k$. 

\begin{equation}
\mathcal{B}^{(k)}=\left\{\boldsymbol{Q}^{(k)},\boldsymbol{c}^{(k)},\boldsymbol{A}^{(k)},\boldsymbol{b}^{(k)},\boldsymbol{G}_0^{(k)},\boldsymbol{h}_0^{(k)},\boldsymbol{G}_i^{(k)},\boldsymbol{h}_i^{(k)},\boldsymbol{a}_i^{(k)},b_i^{(k)},i\in [1, N_{\mathrm{SOC}}]\right\}
\end{equation}

Using a numerical interior-point solver \cite{domahidiECOSSOCPSolver2013}, we obtain the optimal solution of the subproblem, denoted $\boldsymbol{z}^{(k)}_{\star}$. This solution serves simultaneously as the reference trajectory for the next iteration $\boldsymbol{z}_{\mathrm{ref}}^{(k+1)} = \boldsymbol{z}^{(k)}_{\star}$. The SCP forward process is terminated once the improvement in the objective falls below the prescribed threshold given by 
\begin{equation}
    \left|J_2(\boldsymbol{z}^{(k)}_{\star})-J_2(\boldsymbol{z}_{\mathrm{ref}}^{(k)})\right|<\varepsilon_J\left|J_2(\boldsymbol{z}_{\mathrm{ref}}^{(k)})\right|
\end{equation}
where $\varepsilon_J$ is a user-selected convergence tolerance.
A key distinction between the forward process of DSCP and standard SCP lies in the additional information that must be retained for gradient computation. After solving the SOCP, the dual variables associated with inequality constraints must be examined to determine the active sets that enter the backward-pass sensitivity equations. Since numerical solvers return dual values that are not exactly zero even for inactive constraints, Eq.\eqref{equation:active} must be implemented using the thresholds of constraint residuals and dual variables:
\begin{equation}
    \label{equation:active1}
    \begin{aligned}
        \mathcal A_{\mathrm{lin}}
 &\coloneq \left\{ j\left| \, 0 \leq (h_0 - G_0 z_{\star})_j\leq \varepsilon_{\mathrm{res}}, \mu_{j} > \varepsilon_{\mu} \right.\right\}\\
 \mathcal A_{\mathrm{soc}}
 &\coloneq \left\{ i \left| \, 0 \leq g_{i,0} - \left\|\boldsymbol{g}_{i,1}\right\|_2\leq \varepsilon_{\mathrm{res}}, \lambda_{i, 0} \geq \varepsilon_{\lambda}\right.\right\}
    \end{aligned}
\end{equation}
where $\varepsilon_{\mathrm{res}},\varepsilon_{\mu},\varepsilon_{\lambda}$ are empirically selected tolerances informed by the magnitude of residuals and dual variables encountered during the forward iterations. The resulting active sets, along with the extracted quantities $\boldsymbol{\mu}_{+},\boldsymbol{G}_{0+}$, and related matrices required in Eq.\eqref{equation:GS} and Eq.\eqref{equation:KKT1}, are retained in memory for use during backpropagation.  

\subsubsection{Backward Process of Differentiable Sequential Convex Programming}

Assume that the forward SCP process terminates after $K$ iterations. Define a scalar cost function 
\begin{equation}
\mathcal{L}=\mathcal{L}\left(\boldsymbol{z}_{\mathrm{ref}}^{(1)},\dots,\boldsymbol{z}_{\mathrm{ref}}^{(K)},\boldsymbol{\theta}_{\mathrm{SCP}},\boldsymbol{\theta}_{\mathrm{TP}}\right)=\mathcal{L}\left(\boldsymbol{\theta}^{(1)},\dots,\boldsymbol{\theta}^{(K)},\boldsymbol{\theta}^{(K+1)}\right)
\end{equation}
where $\boldsymbol{\theta}^{(K+1)}=\left(\boldsymbol{z}_{\star}^{(K)},\boldsymbol{\theta}_{\mathrm{SCP}},\boldsymbol{\theta}_{\mathrm{TP}}\right)$. The cost function thus depends explicitly on the optimal solution of the final iteration as well as on all intermediate iterations through the sequence of reference trajectories. Since the SCP process is parameterized by $\boldsymbol{\theta}_{\mathrm{SCP}},\boldsymbol{\theta}_{\mathrm{TP}}$, the chain rule gives
\begin{equation}
    \label{equation:d1}
    \begin{aligned}
        \frac{\partial \mathcal{L}}{\partial \boldsymbol{\theta}_{\mathrm{SCP}}}&=\sum_{k=1}^{K+1}\left(\frac{\partial \mathcal{L}}{\partial \boldsymbol{\theta}^{(k)}}\right)^{\top}\frac{\partial \boldsymbol{\theta}^{(k)}}{\partial \boldsymbol{\theta}_{\mathrm{SCP}}}\\
        \frac{\partial \mathcal{L}}{\partial \boldsymbol{\theta}_{\mathrm{TP}}}&=\sum_{k=1}^{K+1}\left(\frac{\partial \mathcal{L}}{\partial \boldsymbol{\theta}^{(k)}}\right)^{\top}\frac{\partial \boldsymbol{\theta}^{(k)}}{\partial \boldsymbol{\theta}_{\mathrm{TP}}}
    \end{aligned}
\end{equation}
The derivative of the cost function with respect to the final optimal solution $\boldsymbol{z}_{\star}^{(K)}$ can be computed directly from the definition of the loss, yielding $\frac{\partial \mathcal{L}}{\partial \boldsymbol{z}^{(K)}}$. For an intermediate iteration $k$, the gradient with respect to the parameter set $\boldsymbol{\theta}^{(\theta)}$ is given by 
\begin{equation}
    \frac{\partial \mathcal{L}}{\partial \boldsymbol{\theta}^{(k)}}=\left(\frac{\partial \mathcal{L}}{\partial \boldsymbol{z}_{\star}^{(k)}}\right)^{\top}\frac{\partial \boldsymbol{z}_{\star}^{(k)}}{\partial \boldsymbol{\theta}^{(k)}},k\in [1,K]
\end{equation}
The Jacobian $\frac{\partial \boldsymbol{z}_{\star}^{(k)}}{\partial \boldsymbol{\theta}^{(k)}}$ can be obtained from the linear sensitivity equation Eq.\eqref{equation:KKT1}. However, computing this Jacobian explicitly would require forming blockwise inverses of $\boldsymbol{H}^{(k)}$ and evaluating explicit analytic derivatives of matrices such as $\boldsymbol{A}^{(k)}$ with respect to $\boldsymbol{\theta}^{(k)}$. To avoid this computational burden, an auxiliary linear system is introduced as 
\begin{equation}
    \label{equation:KKT2}
    \boldsymbol{H}^{(k)\top} \boldsymbol{y}^{(k)} = 
\begin{bmatrix}
\frac{\partial \mathcal{L}}{\partial z_{\star}^{(k)}}^\top & 0 & 0 & 0
\end{bmatrix}^\top
\end{equation}
which has the same coefficient matrix as the forward sensitivity equation and therefore can be solved efficiently using sparse symmetric factorization methods
\begin{equation}    \boldsymbol{P}^{(k)}\boldsymbol{H}^{(k)}\boldsymbol{P}^{(k)\top}=\boldsymbol{L}^{(k)}\boldsymbol{D}^{(k)}\boldsymbol{L}^{(k)^\top}
\end{equation}
Here, $\boldsymbol{P}^{(k)}$ is a permutation matrix obtained using an approximate minimum-degree method, $\boldsymbol{L}^{(k)}$ is lower triangular, and $\boldsymbol{D}^{(k)}$ is a diagonal matrix. A small regularization term is added to the diagonal of $\boldsymbol{H}^{(k)}$ to ensure numerical stability. Since the coefficient matrices arising in SCP are typically large-scale and sparse, exploiting sparse symmetric matrix factorizations to solve the associated linear systems is computationally advantageous and significantly accelerates CPU-based computations. Taken together, Eq.\eqref{equation:KKT1} and Eq.\eqref{equation:KKT2} form a coupled forward-backward propagation system. Eliminating intermediate variables leads to the differential of the loss
\begin{equation}
    \begin{aligned}
    \mathrm{d}\mathcal{L}&=\frac{\partial \mathcal{L}}{\partial z_{\star}^{(k)}}^{\top}\mathrm{d}\boldsymbol{z}_{\star}^{(k)}\\
    &=\begin{bmatrix}\frac{\partial \mathcal{L}}{\partial z_{\star}^{(k)}}^\top &0 &0 &0\end{bmatrix}\begin{bmatrix}\mathrm{d}\boldsymbol{z}_{\star}^{(k)\top}&\text{d}\boldsymbol{\nu}^{\top}&\text{d}\boldsymbol{\mu}_{+}^{\top}&\text{d}\boldsymbol{\lambda}_{0+}^{\top}\end{bmatrix}^\top\\
    &=\begin{bmatrix}\frac{\partial \mathcal{L}}{\partial z_{\star}^{(k)}}^\top &0 &0 &0\end{bmatrix}\boldsymbol{H}^{(k)-1}\boldsymbol{H}^{(k)}\begin{bmatrix}\mathrm{d}\boldsymbol{z}_{\star}^{(k)\top}&\text{d}\boldsymbol{\nu}^{\top}&\text{d}\boldsymbol{\mu}_{+}^{\top}&\text{d}\boldsymbol{\lambda}_{0+}^{\top}\end{bmatrix}^\top\\
    &=\boldsymbol{y}^{(k)\top}\boldsymbol{r}^{(k)}_{\Delta}
    \end{aligned}
\end{equation}
Let the auxiliary vector $\boldsymbol{y}^{(k)\top}$be partitioned as$\left[\boldsymbol{y}_{\boldsymbol{z}}^{(k)\top},\boldsymbol{y}_{\boldsymbol{\nu}}^{(k)\top},\boldsymbol{y}_{\boldsymbol{\mu}}^{(k)\top},\boldsymbol{y}_{\boldsymbol{\lambda}}^{(k)\top}\right]$. Then, the following gradient components can be read directly:
\begin{equation}
    \label{equation:d2}
    \begin{aligned}
    &\frac{\partial \mathcal{L}}{\partial \boldsymbol{A}} = -\left( \boldsymbol{\nu} \boldsymbol{y}_{\boldsymbol{z}}^{\top} + \boldsymbol{y}_{\boldsymbol{\nu}} \boldsymbol{z}^{\top}_{\star} \right),\,\, \frac{\partial \mathcal{L}}{\partial \boldsymbol{G}_i}=-\lambda_{i,0}\left(\frac{\boldsymbol{g}_{i,1}\boldsymbol{y}_{\boldsymbol{z}}^{\top}}{g_{i,0}}+\left(\boldsymbol{I}-\frac{\boldsymbol{g}_{i,1}\boldsymbol{g}_{i,1}^{\top}}{g_{i,0}^2}\right)\frac{\boldsymbol{G}_i\boldsymbol{y}_{\boldsymbol{z}}\boldsymbol{z}^{\top}_{\star}}{g_{i,0}}\right)-y_{s,i}\frac{\boldsymbol{g}_{i,1}\boldsymbol{z}^{\top}_{\star}}{g_{i,0}}
        \\&\frac{\partial \mathcal{L}}{\partial \boldsymbol{b}_0} = \boldsymbol{y}_{\boldsymbol{\nu}},\,\,\frac{\partial \mathcal{L}}{\partial \boldsymbol{h}_{0+}} = \boldsymbol{y}_{\boldsymbol{\mu}},\,\,\frac{\partial\mathcal{L}}{\partial b_i}=y_{s,i},\,\, \frac{\partial \mathcal{L}}{\partial \boldsymbol{h}_i}= -\lambda_{i,0}\left(\boldsymbol{I}-\frac{\boldsymbol{g}_{i,1}\boldsymbol{g}_{i,1}^{\top}}{g_{i,0}^2}\right)\frac{\boldsymbol{G}_i\boldsymbol{y}_{\boldsymbol{z}}}{g_{i,0}}-y_{s,i}\frac{\boldsymbol{g}_{i,1}}{g_{i,0}}
        \\&\frac{\partial \mathcal{L}}{\partial \boldsymbol{G}_{0+}} = \left( \boldsymbol{\mu} \boldsymbol{y}_{\boldsymbol{z}}^{\top}+\boldsymbol{y}_{\boldsymbol{\mu}} \boldsymbol{z}_{\star}^{\top}\right),\,\,\frac{\partial\mathcal{L}}{\partial\boldsymbol{a}_i}=-\lambda_{i,0}\boldsymbol{y}_{\boldsymbol{z}}+y_{s,i}\boldsymbol{z}_{\star},\,\,\frac{\partial \mathcal{L}}{\partial \boldsymbol{Q}}=-\boldsymbol{z}_{\star}\boldsymbol{y}_{\boldsymbol{z}}^{\top},\,\,\frac{\partial \mathcal{L}}{\partial \boldsymbol{c}}=-\boldsymbol{y}_{\boldsymbol{z}}
    \end{aligned}
\end{equation}
These expressions propagate the loss gradient from the primal-dual optimal variables $\boldsymbol{z}_{\star}^{(k)}$ to the collection of constraint and cost coefficients $\mathcal{B}^{(k)}$. Finally, the gradient with respect to the parameters $\boldsymbol{\theta}^{(k)}$ at iteration $k$ is obtained via the chain rule: 
\begin{equation}
    \label{equation:d3}
    \frac{\partial \mathcal{L}}{\partial \boldsymbol{\theta}^{(k)}}=\sum_{\boldsymbol{B}\in \mathcal{B}^{(k)}}\left(\frac{\partial \mathcal{L}}{\partial \boldsymbol{B}}\right)^{\top} \frac{\partial \boldsymbol{B}}{\partial \boldsymbol{\theta}^{(k)}}
\end{equation}
Eq.\eqref{equation:d1}, Eq.\eqref{equation:d2}, and Eq.\eqref{equation:d3} describe the complete backpropagation pathway of the DSCP framework. 

\section{Simulations and Results}
\label{ch4}

To demonstrate the generality and practical utility of the proposed DSCP framework, three numerical experiments are conducted across two representative trajectory planning scenarios: powered powered guidance (PDG) for reusable rockets and entry flight for hypersonic gliding vehicles (HGVs). These experiments target three different categories of parameters, including nonconvex decision variables, hyperparameters of the SCP procedure, and vehicle-design parameters. All simulations are executed in Python 3.8.20. Gradient backpropagation is implemented using PyTorch 1.12.0, and the subproblems are solved using the ECOS solver \cite{domahidiECOSSOCPSolver2013}. As ECOS currently supports only CPU computation, all simulations are performed on the CPU within one core.

\subsection{Optimization of Terminal Time in Powered Descent Guidance}
This subsection applies the DSCP framework to a classical PDG problem \cite{liFreeFinalTimeFuelOptimal2022a}. The objective is to learn and optimize the terminal time, which is a nonconvex and problem-dependent variable that influences fuel optimality.
\subsubsection{Powered Descent Guidance for Reusable Rockets}
The dynamics of the PDG problem is expressed as
\begin{equation}
    \begin{aligned}
        &\dot{\boldsymbol{r}}=\boldsymbol{v}\\
        &\dot{\boldsymbol{v}}=\frac{\boldsymbol{u}_{\mathrm{PDG}}}{m}+\boldsymbol{g}_E\\
        &\dot{m} = -\frac{\left\|\boldsymbol{u}_{\mathrm{PDG}}\right\|_2}{I_{\mathrm{sp}}g_0}
    \end{aligned}
\end{equation}
The state vector is $\boldsymbol{x}_{\mathrm{PDG}}=\left[\boldsymbol{r}^\top,\boldsymbol{v}^\top,m\right]^\top$, where $\boldsymbol{r}$ and $\boldsymbol{v}$ denote the position and velocity vectors, and $m$ is the rocket mass. $ \boldsymbol{u}_{\mathrm{PDG}}$ denotes the thrust of rockets. The parameter $I_{\mathrm{sp}}$ denotes the specific impulse, and $\boldsymbol{g}_E=[g_0, 0, 0]^\top$ is the gravitational acceleration vector, where $g_0$ is the standard gravitational constant. The descent guidance is subject to the following operational constraints:

\begin{equation}
    \begin{aligned}
    &\boldsymbol{x}_{\mathrm{PDG}}(t_0)=\boldsymbol{x}_{\mathrm{PDG},0}\\
    &\boldsymbol{r}(t_{\mathrm{f}})=\boldsymbol{r}_{\mathrm{f}},\boldsymbol{v}(t_{\mathrm{f}})=\boldsymbol{v}_{\mathrm{f}}\\
    &u_{\min}\leq\left\|\boldsymbol{u}_{\mathrm{PDG}}\right\|_2\leq u_{\max}\\
    &\sqrt{r_y^2+r_z^2}-\tan{\beta_{\max}}r_x\leq 0\\
    &\sqrt{u_y^2+u_z^2}-\tan{\eta_{\max}}u_x\leq 0\\
    \end{aligned}
\end{equation}
where $u_{\mathrm{min}}$ and $u_{\mathrm{max}}$ are thrust bounds, $\beta_{\max}$ is the maximum allowable glideslope angle, and $\eta_{\max}$ is the maximum thrust pointing angle. The performance index for fuel-optimal descent is 
\begin{equation}
    J_{\mathrm{PDG}} = -m(t_{\mathrm{f}})    
\end{equation}
Fuel-optimal powered descent is inherently a free-final-time problem. The present experiment evaluates the DSCP framework's ability to learn and optimize such nonconvex timing parameters. The work in \cite{liFreeFinalTimeFuelOptimal2022a} converts the free-final-time problem into multiple fixed-time problems and performs a grid search over the terminal time, followed by training a neural network to approximate the optimal relationship. In contrast, DSCP uses the gradients obtained by backpropagation through the SCP iterations, allowing direct end-to-end training of a neural network without terminal-time enumeration. 

\subsubsection{Verification of Differentiability}
This experiment evaluates the accuracy of the gradients produced by DSCP. A virtual-control-augmented SCP method is used to solve the fixed-time fuel-optimal PDG problem. Let $K_{\mathrm{PDG}}$ be the total number of SCP iterations and $N_{\mathrm{PDG}}$ be the number of discretization intervals. The scalar cost function used for backpropagation is 
\begin{equation}
    \mathcal{L}_{\mathrm{PDG}} = -m^{(K_{\mathrm{PDG}})}[N_{\mathrm{PDG}}]
\end{equation}
By sweeping the terminal time $t_{\mathrm{f}}$, the gradients computed by DSCP are compared against finite-difference gradients of the performance index. This comparison validates the correctness and numerical reliability of the DSCP backward pass. The initialization strategy for the reference trajectory and the construction of SCP subproblems are summarized in Appendix \ref{appendix:PDG}. The problem parameters are set as shown in Table \ref{table:pdg}, and the initial states are given by
\begin{equation}
    \boldsymbol{r}(t_0)=[5000, 500, 500]^\top(\mathrm{m}),\boldsymbol{v}(t_0)=[-150, 30,-30]^\top(\mathrm{m/s}),t_{\mathrm{f}}\in[32.0,34.0](\mathrm{s})
\end{equation}
\begin{table}[h!]
\centering
\caption{Parameters for PDG}
    \label{table:pdg}
\begin{tabular}{c c c c}
\hline
Parameter & Value & Parameter & Value \\
\hline
$u_{\min}$ & $169.0\ \mathrm{kN}$ & $u_{\max}$ & $845.2\ \mathrm{kN}$ \\
$\eta_{\max}$ & $30^\circ$ & $\beta_{\max}$ & $80^\circ$ \\
$I_{\mathrm{sp}}$ & $282.0\ \mathrm{s}$ & $m_0$ & $38000\ \mathrm{kg}$ \\
$N_{\mathrm{PDG}}$ & $50$ & $g_0$ & $9.80655\ \mathrm{m/s^2}$ \\
\hline
\end{tabular}
\end{table}

\begin{figure*}[htbp]
    \centering
    \begin{subfigure}{0.49\textwidth}
        \centering
        \includegraphics[width=\textwidth]{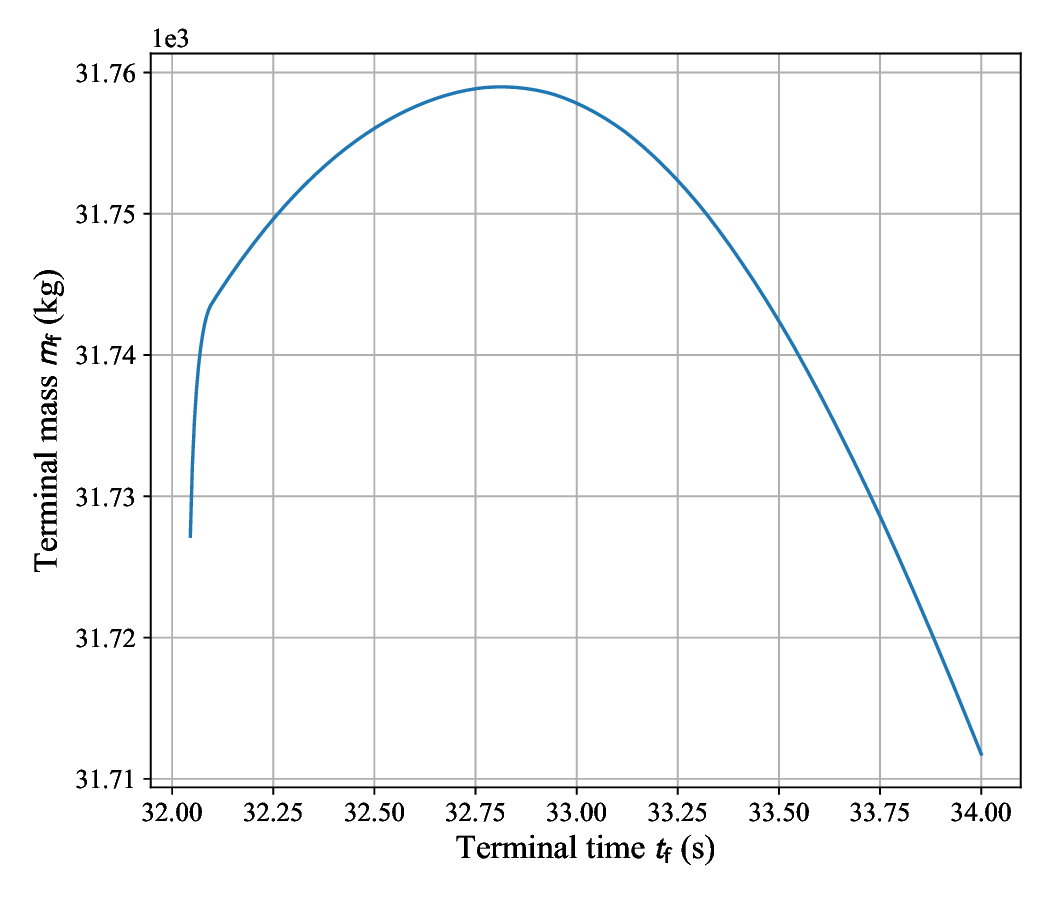}
        \caption{The relationship between maximum remainning mass and terminal time.}
        \label{fig:1-1-1}
    \end{subfigure}
    \hfill
    \begin{subfigure}{0.49\textwidth}
        \centering
        \includegraphics[width=\textwidth]{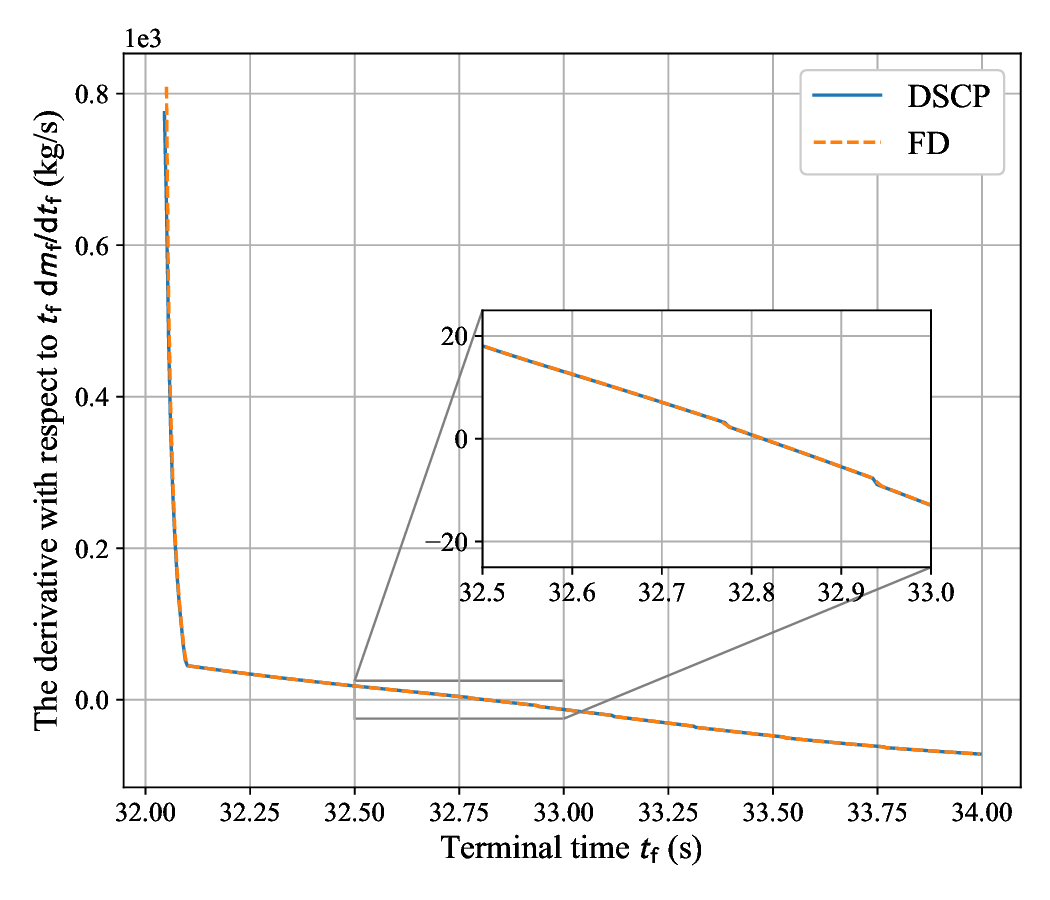}
        \caption{The comparison between gradients obtained by DSCP and finite-difference.}
        \label{fig:1-1-2}
    \end{subfigure}
    \caption{The verification of differentiability of the proposed DSCP.}
    \label{fig:1-1}
\end{figure*}

\begin{figure*}[!ht]
    \centering
    \includegraphics[width=0.95\textwidth]{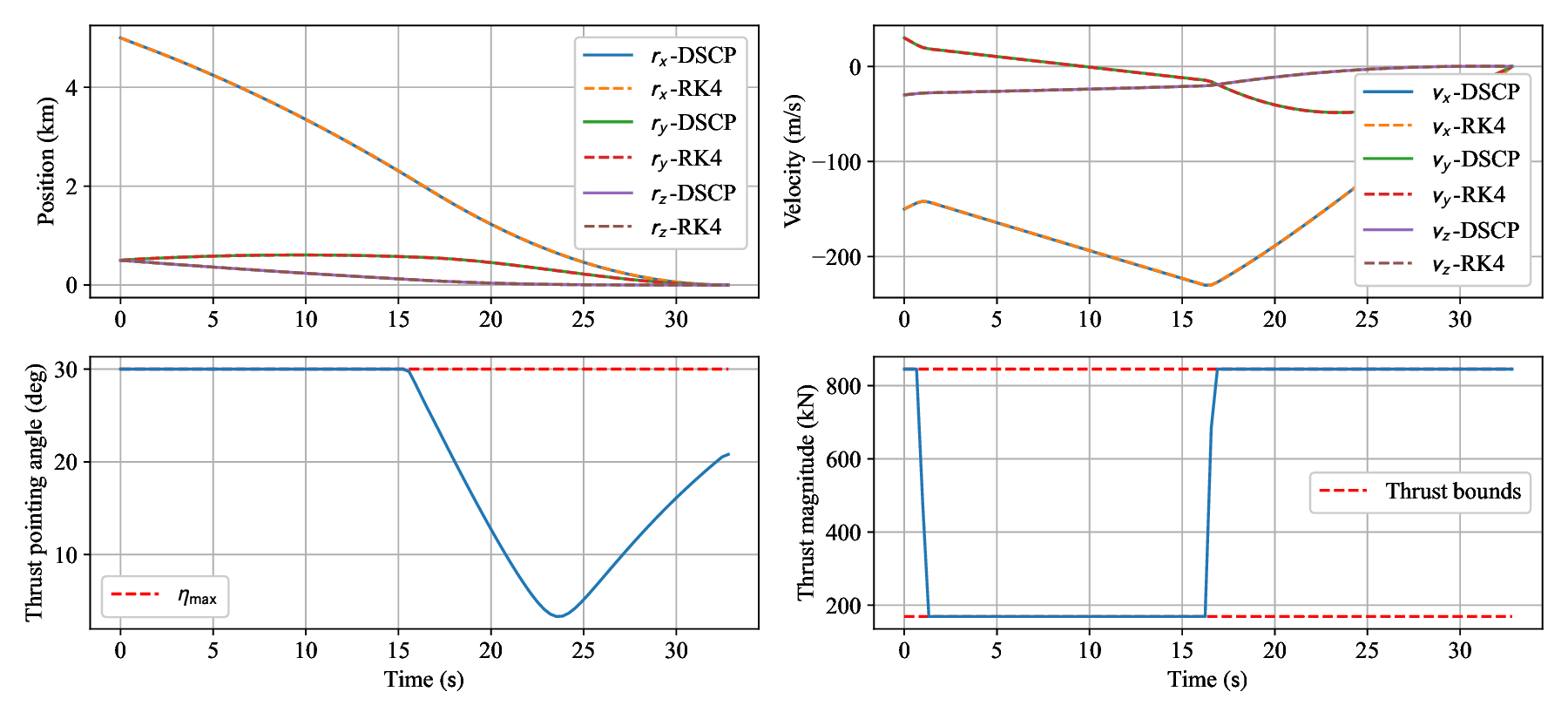}
    \caption{The optimal state profiles computed by DSCP and direct forward simulation at $t_{\mathrm{f}}=32.81$.}
    \label{fig:1-1-3}
\end{figure*}

Simulation results are shown in Fig. \ref{fig:1-1} through Fig. \ref{fig:1-1-3}. Figure \ref{fig:1-1-1} illustrates the relationship between terminal time and maximum remaining mass for a fixed initial condition. The results indicate that no feasible solution exists within the interval $[32.00,32.05]\mathrm{s}$, and solutions in $[32.05,32.10]\mathrm{s}$ exhibit abnormal behavior, where virtual control terms do not converge to zero. This failure occurs because the terminal time is too short for feasible descent. When $t_{\mathrm{f}}>32.10\mathrm{s}$, the SCP iterations converge reliably. With two-decimal precision, the optimal terminal time is $32.81\mathrm{s}$, consistent with the result reported in \cite{liFreeFinalTimeFuelOptimal2022a}. 

Figure \ref{fig:1-1-2} compares gradients computed by DSCP with finite-difference gradients. Except in regions where the subproblem fails to converge, the two gradients match closely, including the zero crossing at $32.81\mathrm{s}$. Figure \ref{fig:1-1-3} compares the optimal solution obtained by DSCP with a direct forward integration of the dynamics by fourth-order Runge-Kutta (RK4). The position, velocity are identical, and thrust magnitude and gimbal angle remain within allowable bounds. These results demonstrate that the SCP formulation used in this work reliably solves the fuel-optimal PDG, and that the proposed DSCP framework provides accurate gradients with respect to problem parameters.

\subsubsection{The Training of Neural Network for Predicting the Terminal Time}
In this experiment, a neural network is trained to predict the terminal time using gradients provided by the DSCP framework. The neural network maps the initial position and velocity of the rockets to the terminal time. The initial states are sampled from
\begin{equation}
    \boldsymbol{r}(t_0)=[5200\pm 200, \pm 600, \pm 600]^\top(\mathrm{m}),\boldsymbol{v}(t_0)=[-135\pm 15, \pm 30,\pm 30]^\top(\mathrm{m/s})
\end{equation}

The network architecture consists of three fully connected layers with 64 neurons per hidden layer and Tanh activation. The final layer uses a sigmoid activation followed by a linear mapping to constrain the output to the interval $[31,37]$s. A plateau-based learning rate decay strategy is adopted, with a decay factor of 0.8. The learning rate is updated every 200 training epochs. In this experiment, the learning rates at epochs 0, 200, and 400 are 0.01, 0.005, and 0.001, respectively. Each epoch contains 256 training samples with a batch size of 32. Since excessively small terminal times lead to infeasible SCP subproblems, the gradient is clipped to $-5.0$ whenever the solver fails. Such failures are frequent early in training, when the neural network's predictions have not yet entered the feasible regime. 

\begin{figure*}[t]
    \centering
    \begin{subfigure}{0.7\textwidth}
        \centering
        \includegraphics[width=\textwidth]{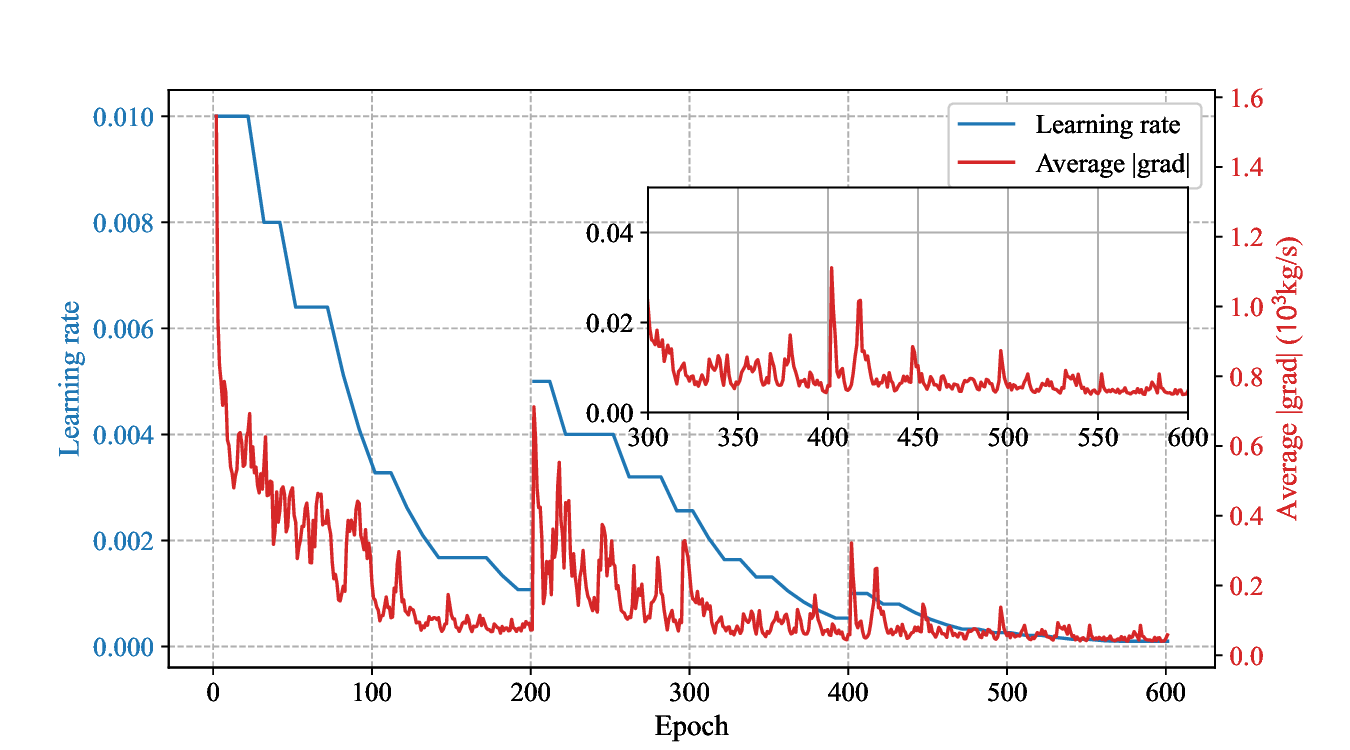}
    \end{subfigure}
    \caption{The training evolution for terminal time prediction using a neural network coupled with DSCP. }
    \label{fig:1_2}
\end{figure*}

\begin{figure*}[!h]
    \centering
    \begin{subfigure}{0.45\textwidth}
        \centering
        \includegraphics[width=\textwidth]{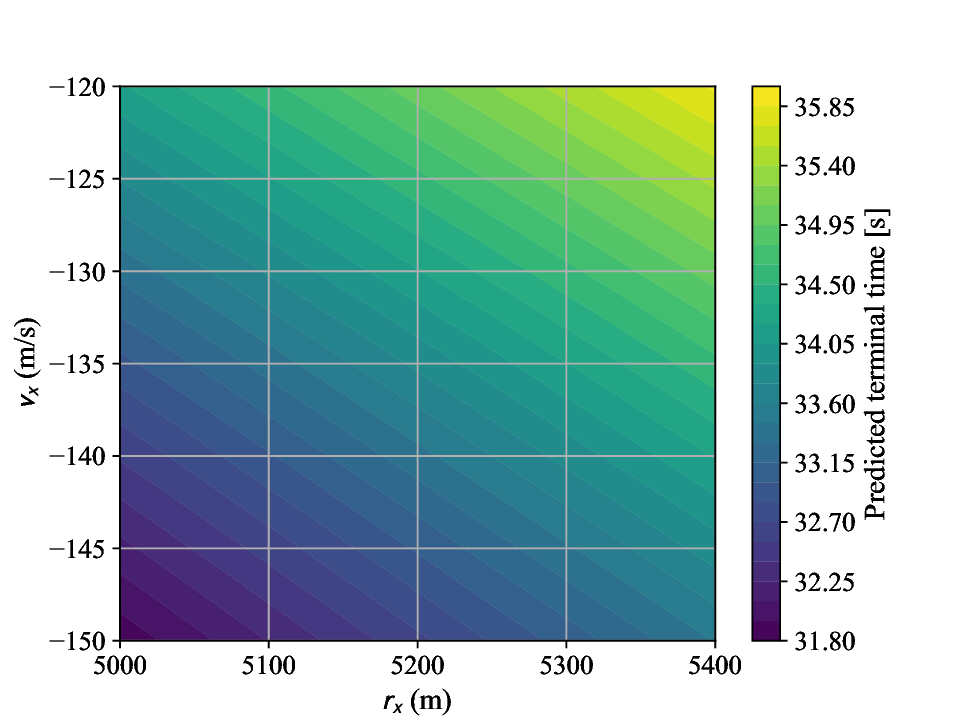}
        \caption{Prediction of the terminal time on the vertical velocity-altitude profile.}
    \end{subfigure}
    \hfill
    \begin{subfigure}{0.45\textwidth}
        \centering
        \includegraphics[width=\textwidth]{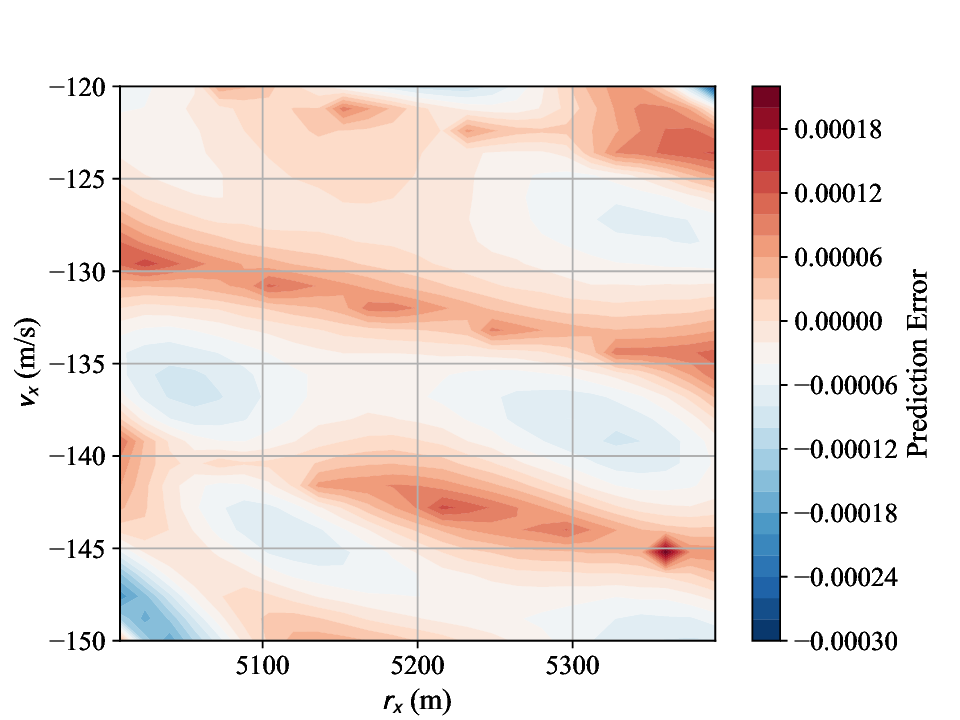}
        \caption{Relative error of the terminal time on the vertical velocity-altitude profile.}
    \end{subfigure}
    \caption{Test performance on the vertical velocity-altitude profile.}
    \label{fig:1_3}
\end{figure*}

\begin{figure*}[!h]
    \centering
    \begin{subfigure}{0.45\textwidth}
        \centering
        \includegraphics[width=\textwidth]{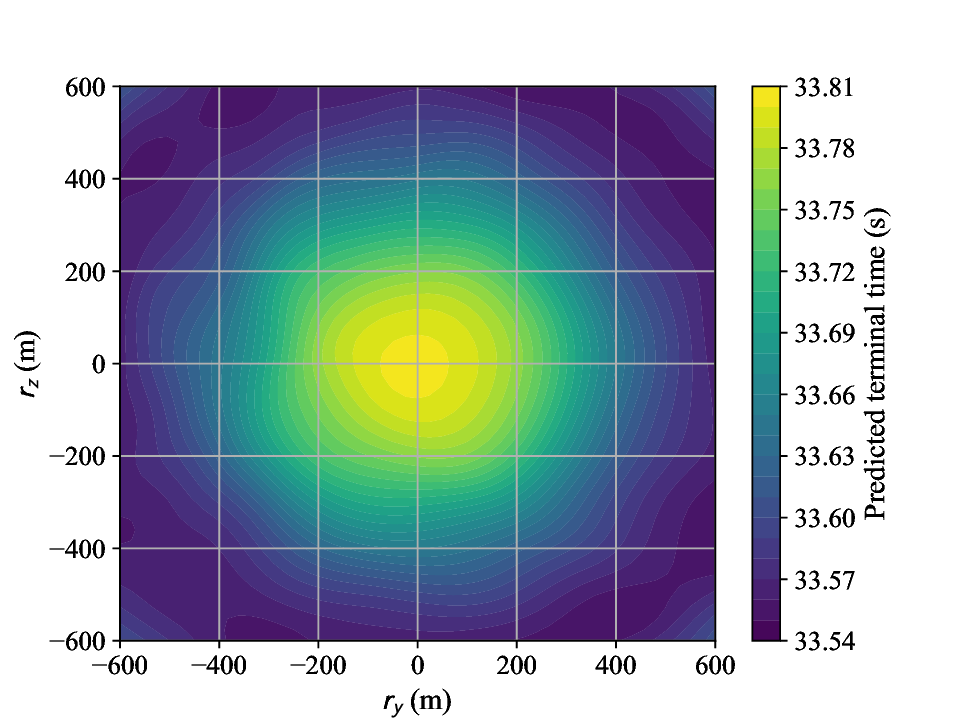}
        \caption{Prediction of the terminal time on the horizontal state profile.}
    \end{subfigure}
    \hfill
    \begin{subfigure}{0.45\textwidth}
        \centering
        \includegraphics[width=\textwidth]{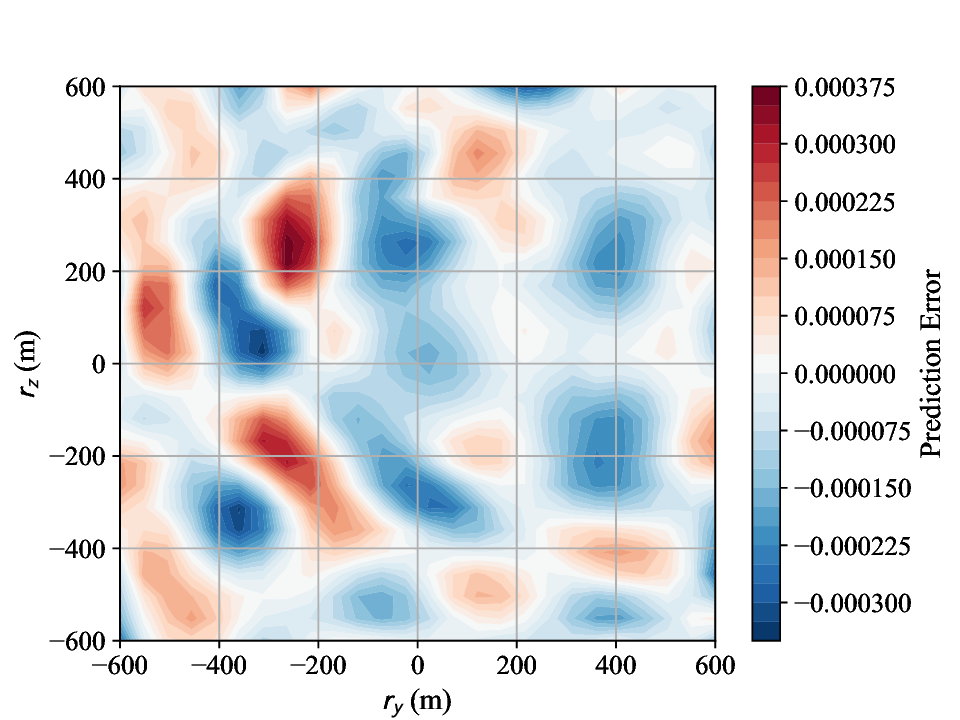}
        \caption{Relative error of the terminal time on the horizontal state profile.}
    \end{subfigure}
    \caption{Test performance on the horizontal state profile.}
    \label{fig:1_4}
\end{figure*}

Figure \ref{fig:1_2} displays the training behavior, including the learning-rate profile and the average absolute value of the output gradient per epoch. As training progresses, gradient magnitudes decrease steadily, indicating improved numerical stability. Test results are presented in Figs.\ref{fig:1_3} and \ref{fig:1_4}. Under a representative vertical velocity-altitude profile with $r_y=r_z=100\mathrm{m},v_y=v_z=0$, the maximum relative prediction error is below $3\times 10^{-4}$. For a horizontal profile with $r_x=5200\mathrm{m},v_x=-135\mathrm{m/s},v_y=v_z=0$, the maximum relative error remains below $4\times 10^{-4}$. These results confirm that the proposed DSCP framework can be seamlessly integrated with neural networks and can effectively train them to learn variables such as the optimal terminal time in trajectory planning problems.

\subsection{Trust-Region Penalty Coefficient Optimization}
\label{ch4-2}

The trust-region penalty coefficients constitute a set of critical design parameters in construction of the convexified subproblems of SCP. This subsection considers the entry trajectory planning of HGVs and investigates the effects of the magnitudes of the trust-region penalty coefficients on the convergence rate and numerical behavior of the trajectory optimization process. In contrast to adaptive trust-region updating strategies used in previous studies, the present investigation focuses on fixed penalty values that remain constant throughout all SCP iterations.


\subsubsection{Entry Trajectory Planning for Hypersonic Gliding Vehicle}
The dynamics of HGV is expressed as \cite{zhangReentryTrajectoryOptimization2025}
\begin{equation}
    \begin{aligned}
        &\dot{r}=V\sin\gamma\\&\dot{\chi}=\frac{V\cos\gamma\sin\psi}{r\cos\phi}\\&\dot{\phi}=\frac{V\cos\gamma\cos\psi}{r}\\
        &\dot{V}=-\frac{D}{m}-g\sin\gamma+r\omega_E^2\cos\phi(\cos\phi\sin\gamma-\sin\phi\cos\gamma\cos\psi)\\
        &\dot{\gamma}=\frac{L\cos\sigma}{mV}-\frac{g\cos\gamma}{V}+\frac{r\omega_E^2\cos\phi}{V}(\cos\phi\cos\gamma+\sin\phi\sin\gamma\cos\psi)+\frac{V\cos\gamma}{r}+2\omega_E\cos\phi\sin\psi\\
        &\dot{\psi}=\frac{L\sin\sigma}{mV\cos\gamma}+\frac{r\omega_{E}^{2}\cos\phi\sin\phi\sin\psi}{V\cos\gamma}+\frac{V\sin\psi\tan\phi\cos\gamma}{r}+2\omega_E(\sin\phi-\cos\phi\tan\gamma\cos\psi)\\&\dot{\alpha}=u_1\\&\dot{\sigma}=u_2
    \end{aligned}
\end{equation}
where the state vector is $\boldsymbol{x}_{\mathrm{HGV}}=[r,\lambda,\phi,V,\gamma,\psi,\alpha,\sigma]^T$ representing the radial distance, longitude, latitude, velocity magnitude, flight-path angle, heading angle, angle of attack, and bank angle. The control input is $\boldsymbol{u}_{\mathrm{HGV}}=[u_1,u_2]^T$ denoting the rate of the angle of attack and bank angle. The planetary rotation rate is denoted by $\omega_E$. The aerodynamic forces acting on the vehicle are given by
\begin{equation}
    \begin{aligned}
        &L=\rho V^2S_{\mathrm{ref}}C_L/2\\
        &D=\rho V^2S_{\mathrm{ref}}C_D/2\\
        &\rho=\rho_0e^{-h/h_\mathrm{s}}
    \end{aligned}
\end{equation}
where $C_L$ and $C_D$ denote lift and drag coefficients, $S_{\mathrm{ref}}$ is the reference area, $\rho$ is the atmospheric density, and $\rho_0, h_s$ are density model constants. The HGV is subject to boundary conditions, dynamic pressure, heat rate, and normal load constraints:
\begin{equation}
\begin{aligned}
        &\boldsymbol{x}_{\mathrm{HGV}}(t_0) = \boldsymbol{x}_{\mathrm{HGV},0}, V(t_\mathrm{f})\geq V_{\mathrm{f,min}}\\
        &\boldsymbol{x}_{\mathrm{HGV}}(t_\mathrm{f}) = \boldsymbol{x}_{\mathrm{HGV,f}} (\mathrm{except\,for\,} V(t_\mathrm{f}))\\
        &\bar p=0.5\rho V^2\leqslant \bar p_{\max}\\
        &\bar q=k_q\rho^{0.5}V^{3.15}\leqslant \bar q_{\max}\\
        &\bar n=\frac{\sqrt{L^2+D^2}}{mg_0}\leqslant \bar n_{\max}
\end{aligned}   
\end{equation}
where $k_q$ is the heat rate constant, and $g_0$ is the gravitational acceleration. The angle and angular-rate constraints are
\begin{equation}
    \label{equation:angleconstraint}
    \begin{aligned}
                &\alpha_{\min}\leq \alpha \leq \alpha_{\max},\sigma_{\min}\leq \sigma \leq \sigma_{\max}\\
        &\dot{\alpha}_{\min}\leq \dot{\alpha} \leq \dot{\alpha}_{\max},\dot{\sigma}_{\min}\leq \dot{\sigma} \leq \dot{\sigma}_{\max}
    \end{aligned}
\end{equation}

The performance index is defined as
\begin{equation}
    J_{\mathrm{HGV}} = \int_{t_0}^{t_{\mathrm{f}}} C_1\gamma^2+C_2\left(\dot{\alpha}^2+\dot{\sigma}^2\right)\mathrm{d}t
\end{equation}
where $C_1 $ and $C_2$ are constants. When linearized and discretized along the reference trajectory, it becomes 
\begin{equation}
    \label{equation:J1HGV}
    \begin{aligned}
    J_{1,\mathrm{HGV}} &= T_{\mathrm{ref}}\sum_{n=0}^{N}\left(2C_1\gamma_{\mathrm{ref}}[n] \delta \gamma[n]+2C_2\dot{\alpha}_{\mathrm{ref}}[n]\delta\dot{\alpha}[n] + 2C_2\dot{\sigma}_{\mathrm{ref}}[n]\delta\dot{\sigma}[n]\right)\\
    &+ \delta T \sum_{n=0}^{N}\left(C_1\gamma^2_{\mathrm{ref}}[n] + C_2\dot{\alpha}^2_{\mathrm{ref}}[n] + C_2\dot{\sigma}^2_{\mathrm{ref}}[n]\right)
    \end{aligned}
\end{equation}
Incorporating the trust-region penalty $J_{\mathrm{trust,HGV}}$, the SCP subproblem in iteration $k$ minimizes 
\begin{equation}
    \label{equation:j2hgv}
    \begin{aligned}
    J_{2,\mathrm{HGV}}&=J_{1,\mathrm{HGV}} + J_{\mathrm{trust,HGV}}\\
    & =J_{1,\mathrm{HGV}} + \sum_{n=0}^{N}\left(\omega_{\gamma} \delta \gamma^2[n] + \omega_{u}\left(\delta \dot {\alpha}^2[n] + \delta \dot{\sigma}^2[n]\right)\right)
    \end{aligned}
\end{equation}
where $\omega_{\gamma}$ and $\omega_{u}$ denote the trust-region penalty coefficients. Let $J_{1,\mathrm{HGV}}^{(k)}$ denote the objective value at iteration $k$ and define $K_{\mathrm{HGV}}$ as the total number of iterations before convergence or termination. To quantitatively assess the convergence behavior under different trust-region parameter values, define the evaluation metric
\begin{equation}
    \mathcal{L}_{\mathrm{HGV}}=\sum_{k=2}^{K_{\mathrm{HGV}}} |J_{1,\mathrm{HGV}}^{(k)}-J_{1,\mathrm{HGV}}^{(K_{\mathrm{HGV}})}|
\end{equation}
which indicates that a smaller $\mathcal{L}_{\mathrm{HGV}}$ indicates faster and smoother convergence. Excessively large trust-region penalty coefficients significantly slow the updates of the state and control variables, thereby degrading the convergence performance. Consequently, unreasonable penalty values are excluded from the subsequent simulations. Additional details of the problem setup are provided in Appendix~\ref{appendix:HGV}.

\subsubsection{Optimization of Trust-Region Penalty Coefficients}
The aerodynamic coefficients of HGVs are expressed as
\begin{equation}
\begin{aligned}
C_L(\alpha,Ma) &= C_{L11}\alpha + C_{L21}\alpha^2 + C_{L22}\alpha Ma, \\
C_D(\alpha,Ma) &= C_{D01} + C_{D11}Ma + C_{D21}\alpha^2 + C_{D22}Ma^2,
\end{aligned}
\end{equation}
To fully specify the HGV entry problem, Table~\ref{tab:aero_coeff} lists the aerodynamic coefficients, physical constants, and initial boundary conditions. For clarity, the allowable bounds on the angle of attack, bank angle, and their corresponding rate limits are also included in the table.

\begin{table}[htbp]
\centering
\caption{Parameters for Entry Trajectory Planning}
\label{tab:aero_coeff}
\begin{tabular}{cccccccc}
\hline
Parameter & Value & Parameter & Value & Parameter & Value & Parameter & Value \\
\hline
$C_{L11}$ & 2.497      & $r_0$        & $6449.1~\mathrm{km}$  & $r_{\mathrm{f}}$        & $6401.1~\mathrm{km}$ &$\alpha_{\mathrm{min}}$&$0.0^\circ$\\
$C_{L21}$ & 1.477      & $\chi_0$  & $0.0^\circ$          & $\chi_{\mathrm{f}}$ & $120.0^\circ$ &$\alpha_{\mathrm{max}}$&$40.0^\circ$\\
$C_{L22}$ & $-0.03731$ & $\phi_0$     & $0.0^\circ$          & $\phi_{\mathrm{f}}$    & $12.0^\circ$ &$\sigma_{\mathrm{min}}$&$-60.0^\circ$\\
$C_{D01}$ & 0.2298     & $V_0$        & $6.7~\mathrm{km/s}$  & $V_{\mathrm{f,min}}$       & $1.0~\mathrm{km/s}$ &$\sigma_{\mathrm{max}}$&$60.0^\circ$\\
$C_{D11}$ & $-0.02432$ & $\gamma_0$   & $0.0^\circ$          & $\gamma_{\mathrm{f}}$  & $0.0^\circ$ &$\dot{\alpha}_{\mathrm{min}}$&$-5.0^\circ\mathrm{/s}$\\
$C_{D21}$ & 2.36       & $\psi_0$     & $90.0^\circ$         & $\psi_{\mathrm{f}}$    & $90.0^\circ$ &$\dot{\alpha}_{\mathrm{max}}$&$5.0^\circ\mathrm{/s}$\\
$C_{D22}$ & 0.0007079  & $\alpha_0$   & $30.0^\circ$         & $t_{\mathrm{f,min}}$  & $1800\mathrm{s}$ &$\dot{\sigma}_{\mathrm{min}}$&$-5.0^\circ\mathrm{/s}$\\
$m_{\mathrm{HGV}}$  & $907.186~\text{kg}$        & $\sigma_0$   & $0.0^\circ$          & $t_{\mathrm{f,min}}$  & $3000\mathrm{s}$ &$\dot{\sigma}_{\mathrm{max}}$&$5.0^\circ\mathrm{/s}$\\
\hline
\end{tabular}
\end{table}
The common path constraints on dynamic pressure, heat rate, and load factor are imposed as $\bar{p}_{\max}=1.1\times10^5~\text{N/m}^2,\bar{q}_{\max}=2.5\times10^6~\text{W/m}^2,\bar{n}_{\max}=2.0$ with reference area $0.48~\mathrm{m^2}$ and heat rate constant $9.437\times10^{-5}$.For numerical stability and scaling purposes, the normalization parameters adopted in the simulations are set as $r_{\mathrm{scale}} = 1.0\times10^5~\text{m},v_{\mathrm{scale}} = 1.0\times10^4~\text{m/s},m_{\mathrm{scale}} = 907.186~\text{kg}$. The initial reference trajectory is generated via forward numerical integration under a prescribed control profile given in Appendix~\ref{appendix:HGV}. The total number of discretization intervals is chosen as $N_{\mathrm{HGV}}=100$, and the SCP convergence tolerance is prescribed as $\varepsilon = 0.001$. Unless otherwise specified, the unified simulation settings are adopted for all numerical cases in this subsection.  

The parameters applied in the performance index Eq.\eqref{equation:J1HGV} are selected as $C_1 = 1.0, C_2 = 1.1$. The trust-region coefficients are sampled from
\begin{equation}
    \omega_{\gamma}, \omega_{u}\in [10^{-1},10^{2}]
\end{equation}

\begin{figure*}[t]
\centering
        \includegraphics[width=\textwidth]{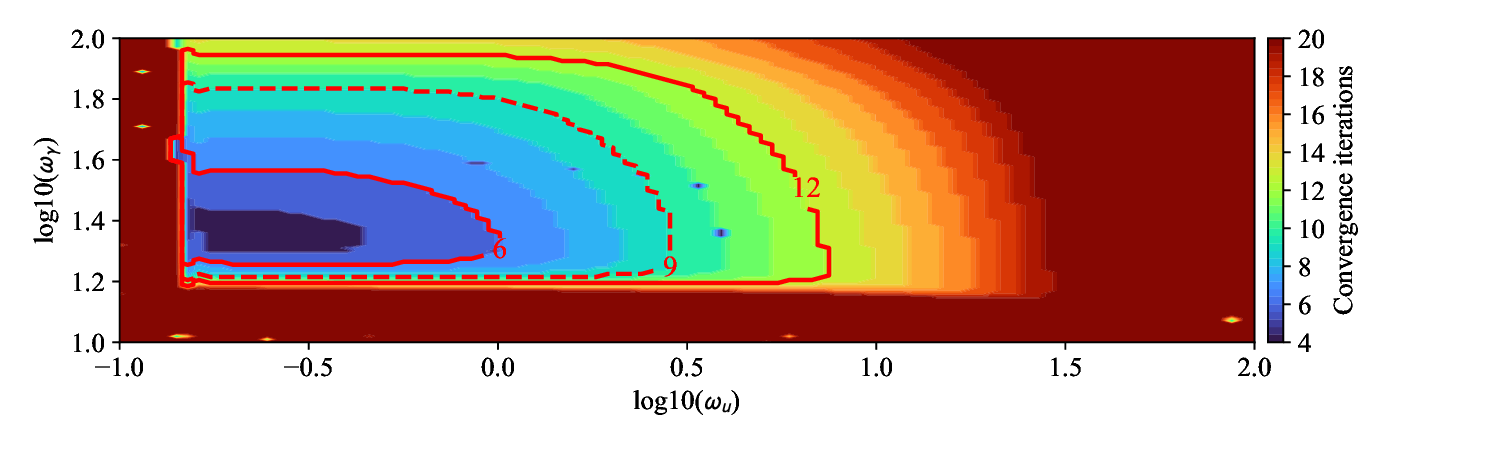}
        \caption{Convergence iterations on the $\omega_{\gamma}-\omega_{u}$ profile.}
        \label{fig:211}
\end{figure*}

\begin{figure*}[t]
\centering
        \includegraphics[width=\textwidth]{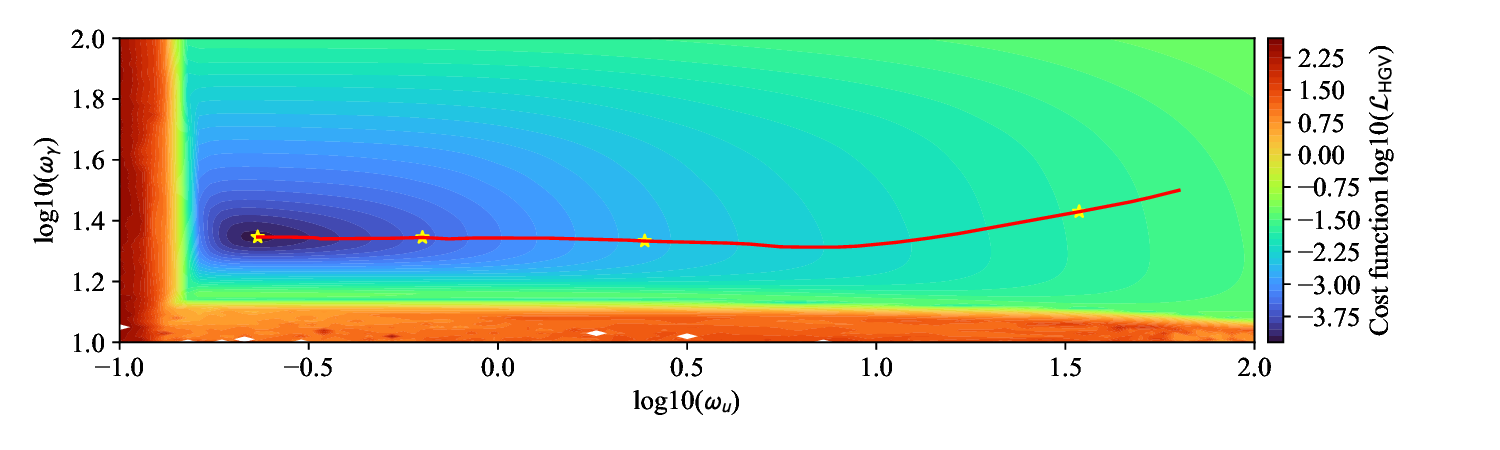}
        \caption{Cost function values on the $\omega_{\gamma}-\omega_{u}$ profile.}
        \label{fig:212}
\end{figure*}

\begin{figure*}[h]
\centering
        \includegraphics[width=\textwidth]{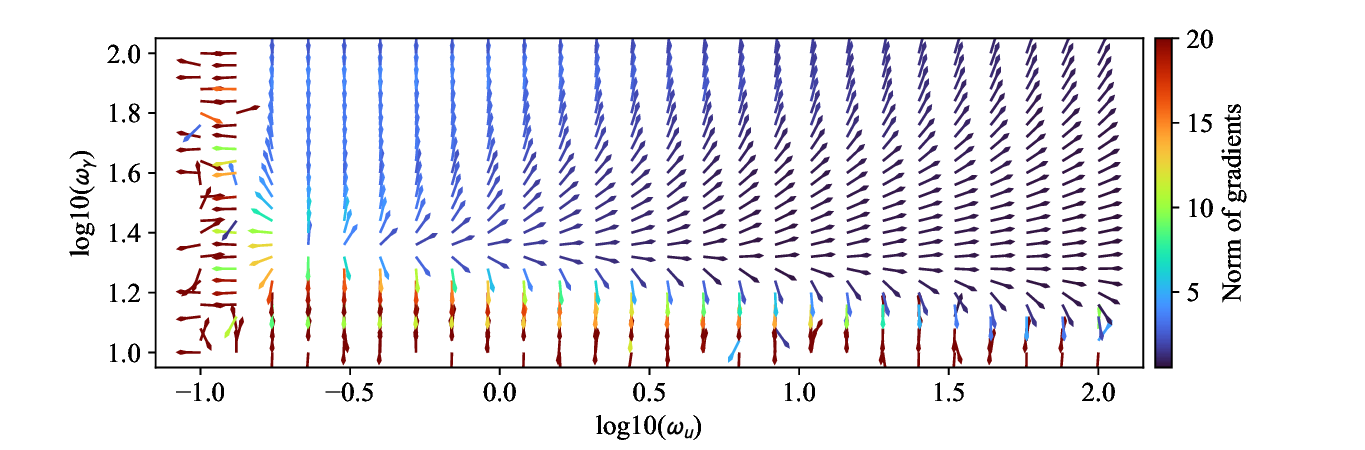}
        \caption{Normalized gradient field with respect to $\omega_{\gamma}$ and $\omega_{u}$}
        \label{fig:213}
\end{figure*}

Simulation results for different trust-region penalties are shown in Figs.~\ref{fig:211} and \ref{fig:212}. Only the region $\omega_{\gamma} \in [10^{1},\,10^{2}], \omega_{u} \in [10^{-1},\,10^{2}]$ is displayed, because smaller $\omega_{\gamma}$ values fail to provide sufficient regularization for the flight-path angle, causing SCP to diverge within the maximum iteration limit (20 iterations). Accordingly, the displayed region represents the primary domain in which SCP approaches or achieves convergence. 

Figure~\ref{fig:211} illustrates the distribution of SCP convergence iterations. A well-defined basin of fast convergence is observed around $\omega_{\gamma}\in[10^{1.3},\,10^{1.5}], \omega_{u}\in[10^{-0.7},\,10^{0.4}]$. Within this region, the trust-region penalties effectively constrain the update magnitudes of state and control variables, thereby preserving the validity of linearization and accelerating convergence. Outside this basin, the number of iterations increases significantly, highlighting the sensitivity of SCP to unbalanced or improperly scaled trust-region penalty parameters.

Figure~\ref{fig:212} presents the corresponding distribution of the cost function. The valley of minimum cost closely aligns with the region of rapid convergence, confirming that the proposed convergence metric provide a meaningful characterization of SCP performance. The red trajectory illustrates the parameter optimization update direction of the algorithm from a representative initial point. It exhibits a continuous descending trend along the cost contours, indicating that the cost function with respect to the penalty parameters possesses a smooth and differentiable structure. This observation further validates that the proposed DSCP framework is amenable to effective optimization in the parameter space. 
Figure~\ref{fig:213} depicts the gradient field obtained by DSCP backpropagation, quantified by $\left\|\nabla_{\omega_{\gamma},\,\omega_{u}} \mathcal{L}_{\mathrm{HGV}}\right\|$. The gradient field exhibits a clear outward divergence from the minimum-cost region, reflecting the structure of the underlying cost landscape and confirming that DSCP accurately captures the sensitivity of the objective with respect to the trust-region penalty parameters. In contrast, for SCP configurations that fail to converge, the gradients follow the expected outward trend only in a global sense while exhibiting localized irregularities in certain regions. These deviations indicate that effective parameter optimization under DSCP still benefits from appropriately chosen initial parameter values. 

\begin{figure*}[!h]
    \centering
    \begin{subfigure}{0.49\textwidth}
        \centering
        \includegraphics[width=\textwidth]{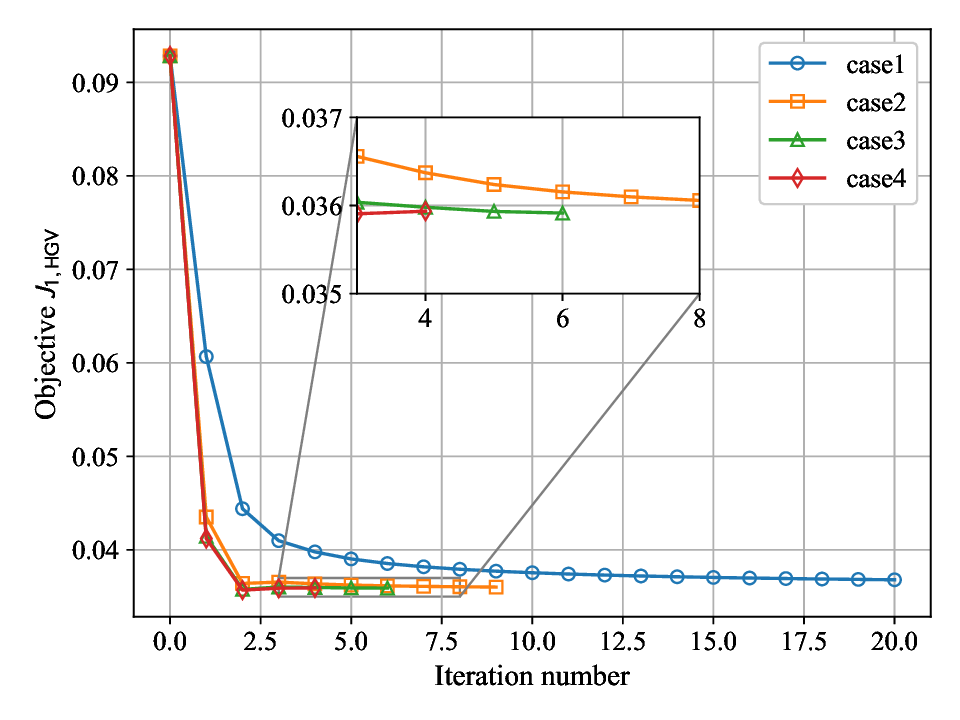}
        \caption{Performance index over DSCP iterations.}
        \label{fig:221}
    \end{subfigure}
    \hfill
    \begin{subfigure}{0.49\textwidth}
        \centering
        \includegraphics[width=\textwidth]{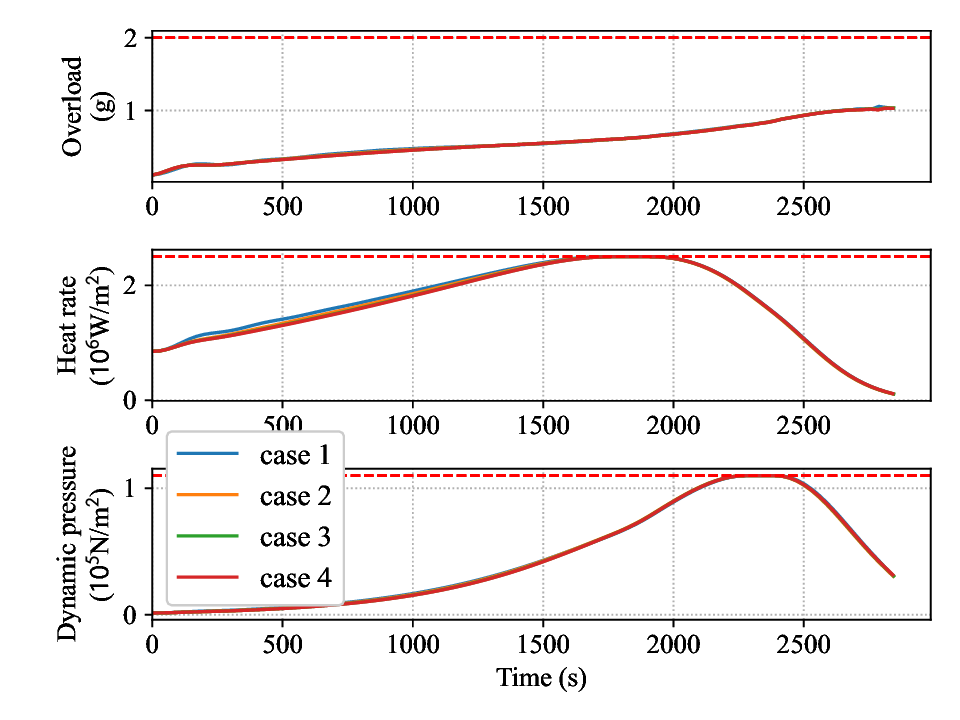}
        \caption{Normal load, heat rate and dynamic pressure.}
    \end{subfigure}
    \caption{Evolution of performance index and constraint satisfaction under different trust-region penalties.}
    \label{fig:2_2}
\end{figure*}

\begin{figure*}[!h]
    \centering
    \begin{subfigure}{0.48\textwidth}
        \centering
        \includegraphics[width=\textwidth]{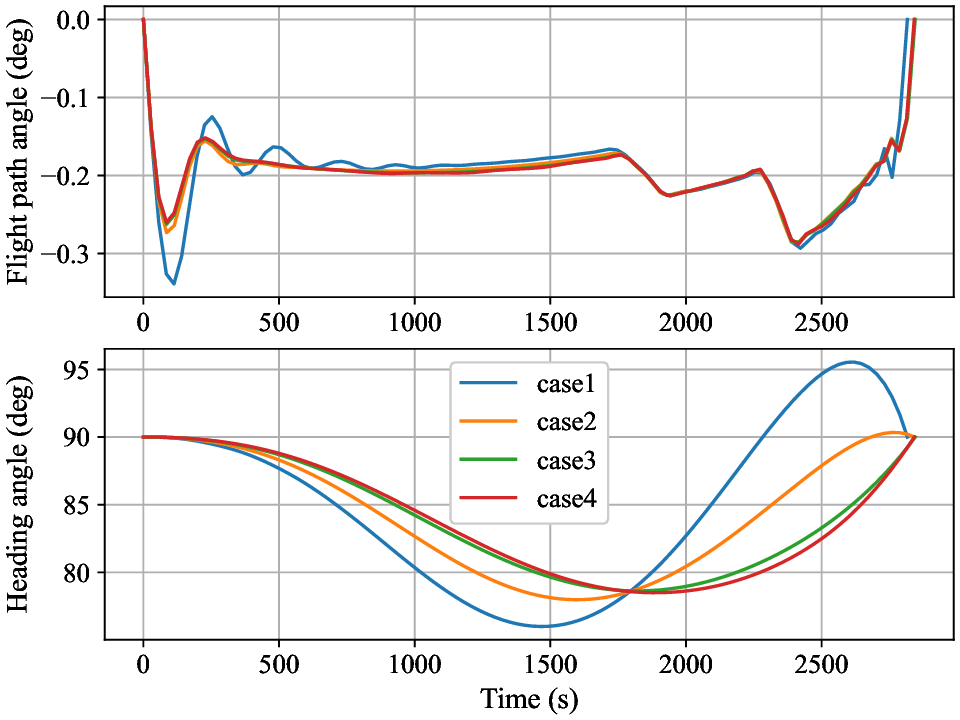}
        \caption{Flight-path angle and heading angle profiles.}
    \end{subfigure}
    \hfill
    \begin{subfigure}{0.48\textwidth}
        \centering
        \includegraphics[width=\textwidth]{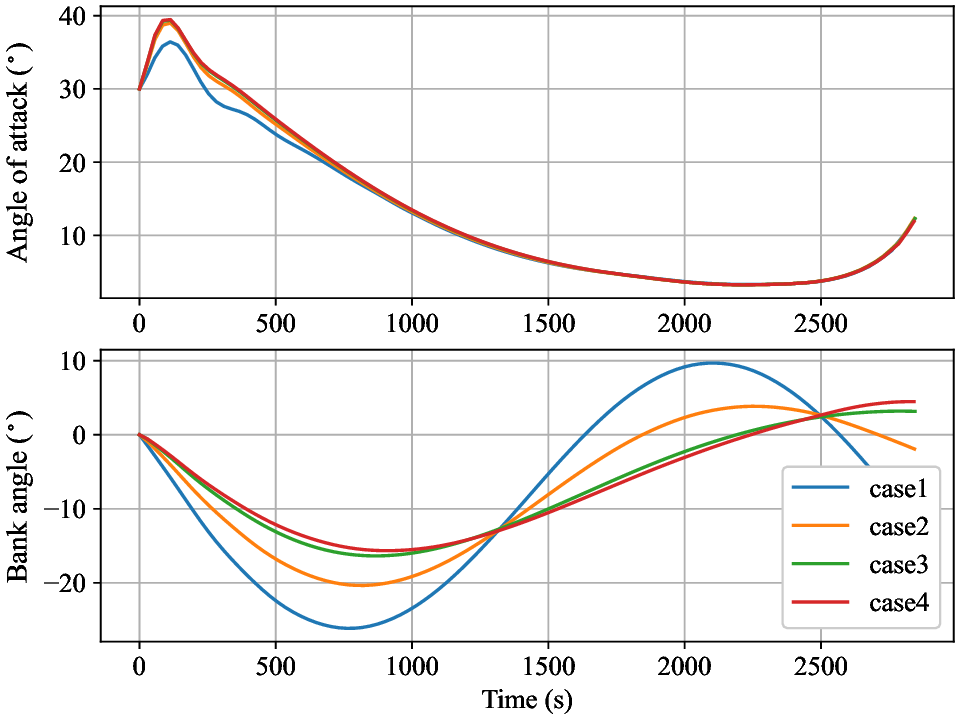}
        \caption{Angle of attack and bank angle profiles.}
    \end{subfigure}
    \caption{State profiles under different trust-region penalties.}
    \label{fig:2_3}
\end{figure*}

\begin{figure*}[!h]
    \centering
    \begin{subfigure}{0.48\textwidth}
        \centering
        \includegraphics[width=\textwidth]{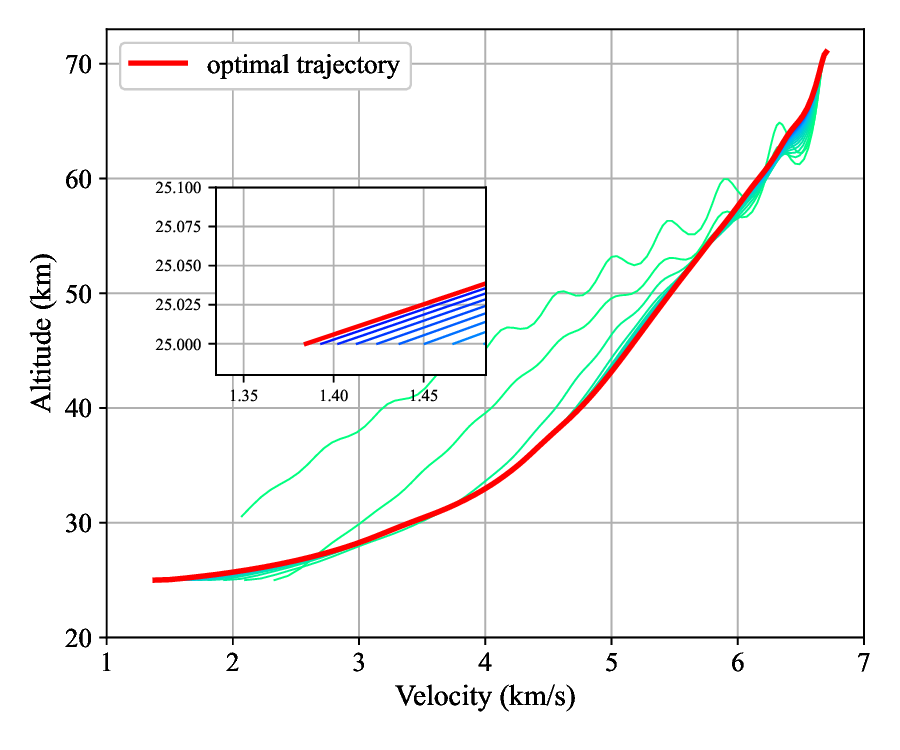}
        \caption{Case1: $\log_{10}(\omega_\gamma)=1.429,\log_{10}(\omega_{u})=1.536$.}
    \end{subfigure}
    \hfill
    \begin{subfigure}{0.48\textwidth}
        \centering
        \includegraphics[width=\textwidth]{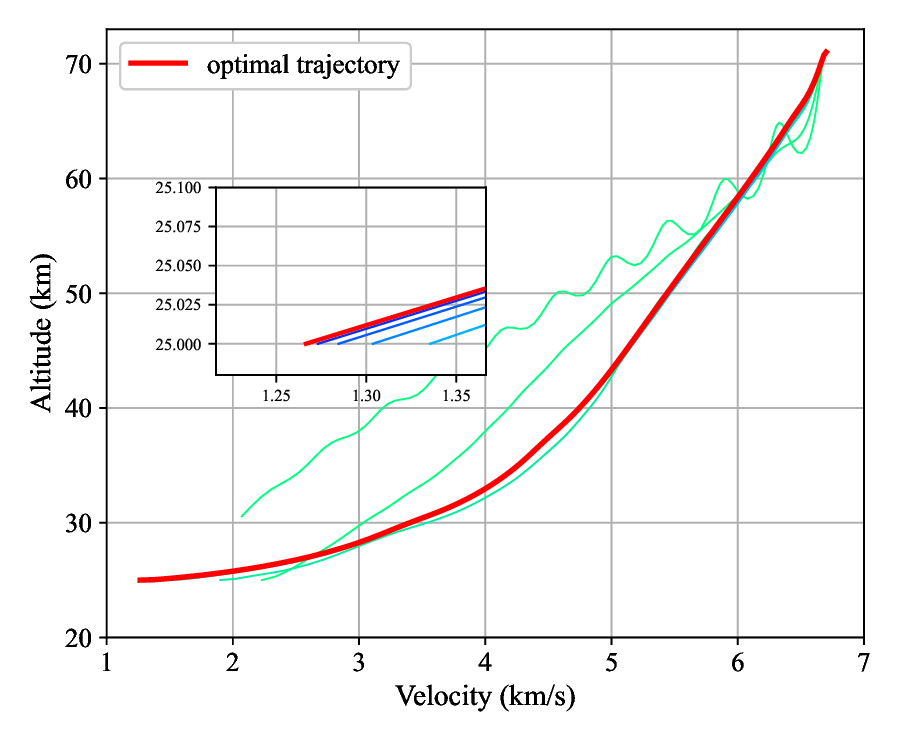}
        \caption{Case2: $\log_{10}(\omega_\gamma)=1.333,\log_{10}(\omega_{u})=0.388$.}
    \end{subfigure}
    \begin{subfigure}{0.48\textwidth}
        \centering
        \includegraphics[width=\textwidth]{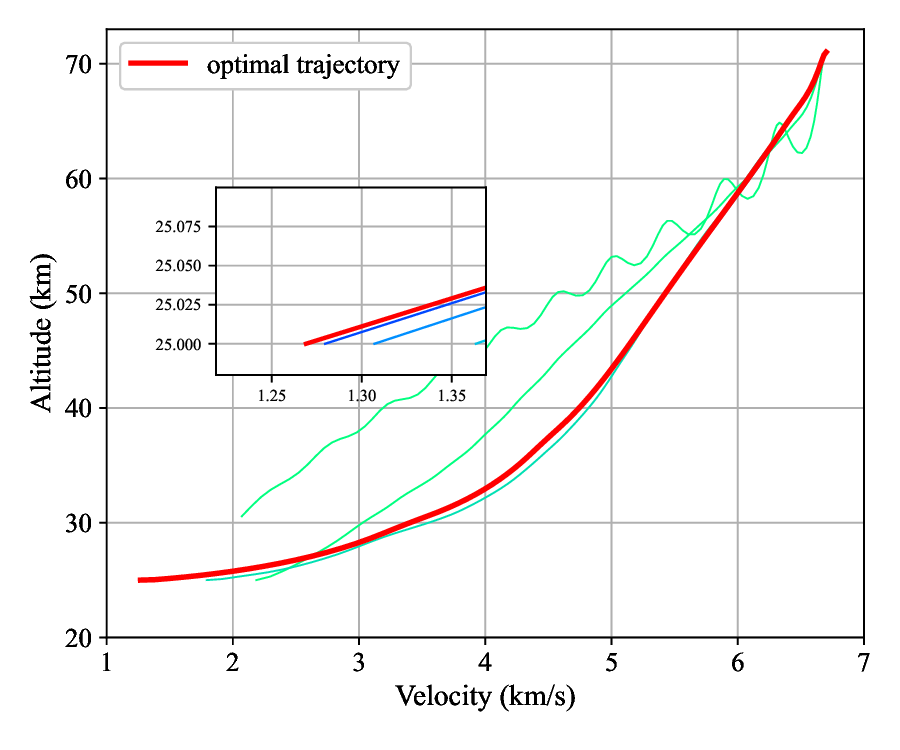}
        \caption{Case3: $\log_{10}(\omega_\gamma)=1.344,\log_{10}(\omega_{u})=-0.199$.}
    \end{subfigure}
    \hfill
    \begin{subfigure}{0.48\textwidth}
        \centering
        \includegraphics[width=\textwidth]{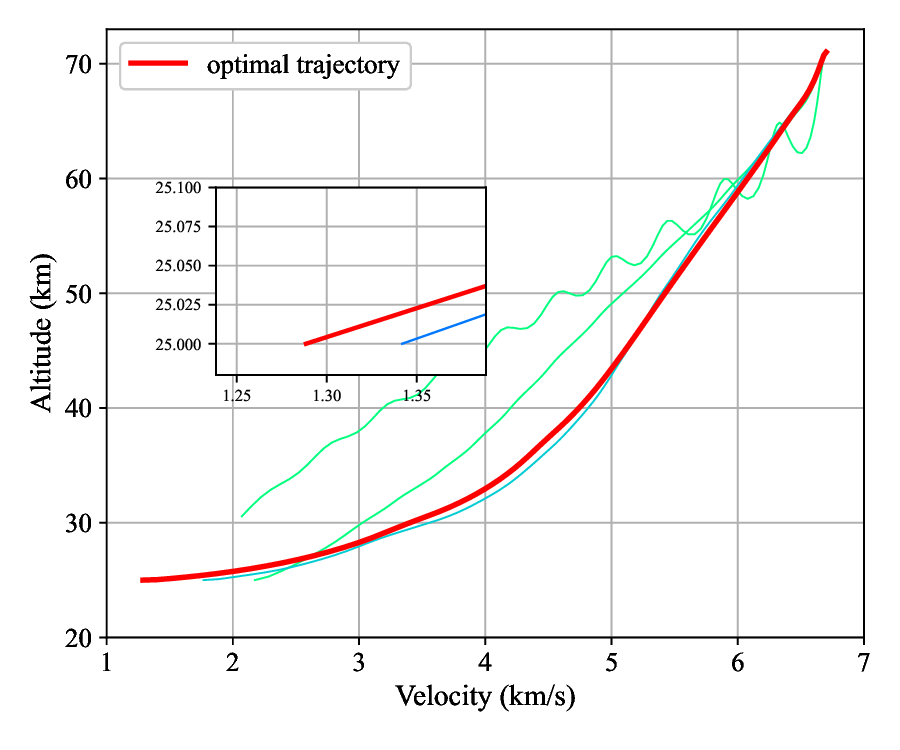}
        \caption{Case4: $\log_{10}(\omega_\gamma)=1.346,\log_{10}(\omega_{u})=-0.635$.}
    \end{subfigure}
    \caption{Altitude-velocity profiles over DSCP iterations under different trust-region penalties.}
    \label{fig:2_4}
\end{figure*}

To further demonstrate the influence of trust-region parameters on the optimization process, four representative cases corresponding to nonconvergence, 9-step convergence, 6-step convergence, and 4-step convergence are selected from the descent curve in Fig.~\ref{fig:212}. The evolution of performance index and constraint satisfaction is illustrated in Fig.~\ref{fig:2_2}. All cases satisfy dynamic pressure, heat rate, and normal load constraints, while cases with faster convergence consistently achieve lower final cost values. Figure~\ref{fig:2_3} displays the evolution of flight-path angle, heading angle, angle of attack, and bank angle. The trajectories differ significantly, and cases with slow convergence tend to evolve suboptimal regions. Figure~\ref{fig:2_4} shows the altitude-velocity profiles for each case. Excessively large trust-region penalties slow down the update of the altitude-velocity profile, thereby degrading the overall convergence rate.

In summary, the simulation results demonstrate the following: (1) the trust-region penalty coefficients exert a significant influence on the convergence behavior and feasibility of SCP; (2) the resulting cost distribution a exhibits smooth and differentiable structure, which enables effective parameter learning; and (3) DSCP accurately captures the relevant gradient information and reflects the sensitivity of the cost function with respect to the penalty parameters. These findings highlight the importance of appropriate trust-region parameter selection and further validate the effectiveness of DSCP for automatic parameter tuning.


\subsection{Geometric Parameters Optimization of HGVs}

This simulation employs the same entry-flight model introduced in the previous subsection and investigates how the geometric properties of HGVs influence the trajectory planning performance. In this simulation, the cost function simplifies to

\begin{equation}
    \mathcal{L}_{\mathrm{HGV}}=J_{\mathrm{1,HGV}}
\end{equation}
The surface-to-mass ratio is a key design parameter for HGVs. Since the aerodynamic coefficients can be further tailored through airframe shaping, the present study neglects the influence of the surface-to-mass ratio on the aerodynamic coefficients themselves. Based on the nominal value $ C_{S/m,0}$ reported in \cite{zhangReentryTrajectoryOptimization2025}, the surface-to-mass ratio is parameterized as
\begin{equation}
    C_{S/m}=\omega_{S/m}C_{S/m,0},\omega_{S/m}\in [0.75, 1.25]
\end{equation}
where $\omega_{S/m}$ denotes the scaling factor. The aim of this experiment is to identify the optimal value of $\omega_{S/m}$ within the interval $[0.75,1.25]$, thereby providing insight into the geometric design characteristics of the vehicles. To further assess the effectiveness of the optimized trust-region coefficients identified, four representative cases of trust-region parameters, which are the same as Sec.~\ref{ch4-2}, are employed for the optimization of $ \omega_{S/m}$. All other simulation settings are kept identical to those used in Sec.~\ref{ch4-2}.


\begin{table}[!t]
\centering
\caption{Statistical performance of surface-to-mass ratio optimization.}
\label{tab:performance}
\begin{tabular}{cccc}
\hline
Case & Optimal parameter  & Optimal value& Iteration number\\
\hline
$N_{\mathrm{HGV}}=100$&&&\\
\hline
case1 & 0.94      & $0.0368$        & $18.92\pm2.49$\\
case2 & 0.93      & $0.0359$        & $7.86\pm2.91$\\
case3 & 0.93      & $0.0358$        & $5.34\pm0.22$\\
case4 & 0.93      & $0.0358$        & $5.13\pm0.73$\\
\hline
$N_{\mathrm{HGV}}=200$&&&\\
\hline
case1 & 1.00      & $0.0373$        & $20.00\pm0.00$\\
case2 & 0.92      & $0.0358$        & $11.59\pm0.40$\\
case3 & 0.91      & $0.0356$        & $8.04\pm0.47$\\
case4 & 0.91      & $0.0356$        & $7.47\pm0.25$\\
\hline
\end{tabular}
\end{table}

\begin{figure*}[t]
    \centering
    \begin{subfigure}{0.49\textwidth}
        \centering
        \includegraphics[width=\textwidth]{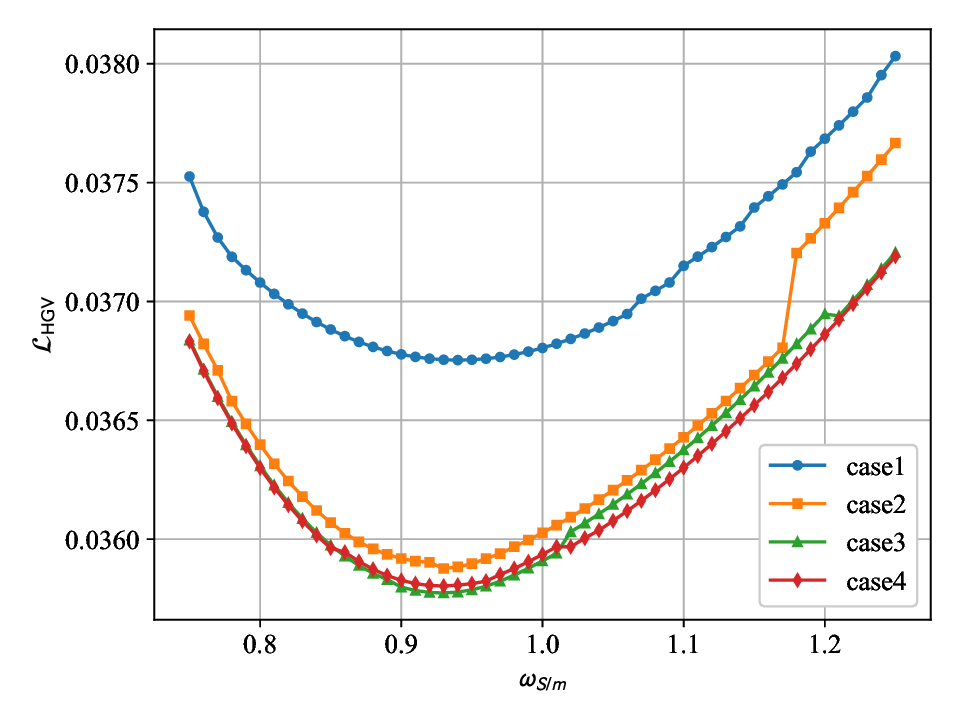}
        \caption{Value of the cost function for different $\omega_{S/m}$.}
        \label{fig:311}
    \end{subfigure}
    \hfill
    \begin{subfigure}{0.49\textwidth}
        \centering
        \includegraphics[width=\textwidth]{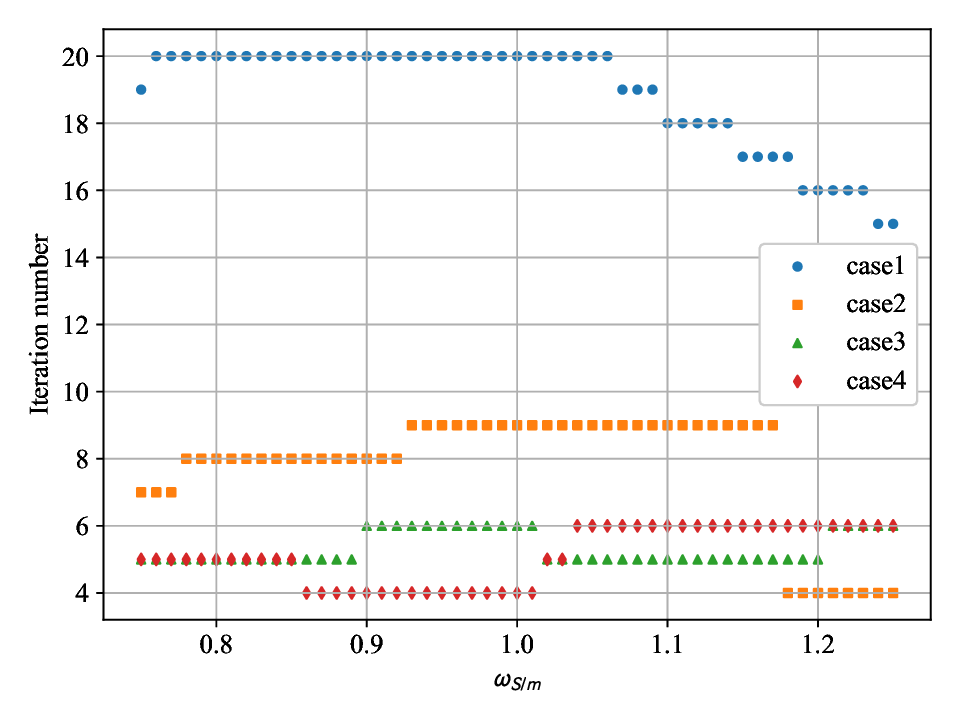}
        \caption{Number of iterations for different $\omega_{S/m}$.}
        \label{fig:312}
    \end{subfigure}
    \caption{Optimization behavior with respect to the surface-to-mass ratio $\omega_{S/m}$ with $N_{\mathrm{HGV}}=100$}
    \label{fig:31}
\end{figure*}

\begin{figure*}[t]
    \centering
    \begin{subfigure}{0.49\textwidth}
        \centering
        \includegraphics[width=\textwidth]{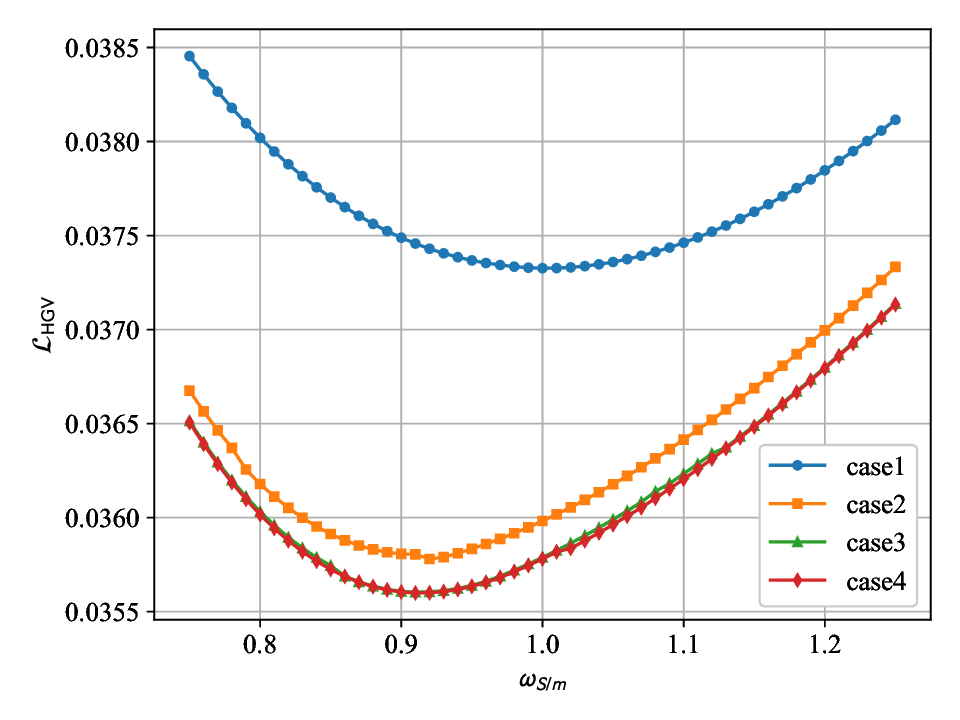}
        \caption{Value of the cost function for different $\omega_{S/m}$.}
        \label{fig:321}
    \end{subfigure}
    \hfill
    \begin{subfigure}{0.49\textwidth}
        \centering
        \includegraphics[width=\textwidth]{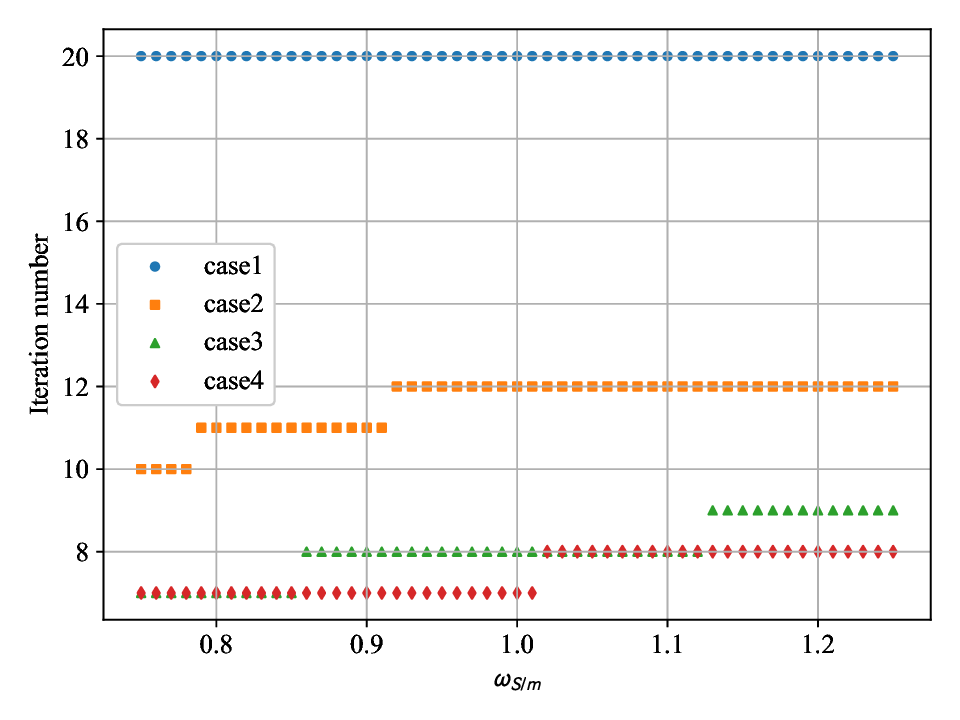}
        \caption{Number of iterations for different $\omega_{S/m}$.}
        \label{fig:322}
    \end{subfigure}
    \caption{Optimization behavior with respect to the surface-to-mass ratio $\omega_{S/m}$ with $N_{\mathrm{HGV}}=200$.}
    \label{fig:32}
\end{figure*}

\begin{figure*}[!h]
    \centering
    \begin{subfigure}{0.49\textwidth}
        \centering
        \includegraphics[width=\textwidth]{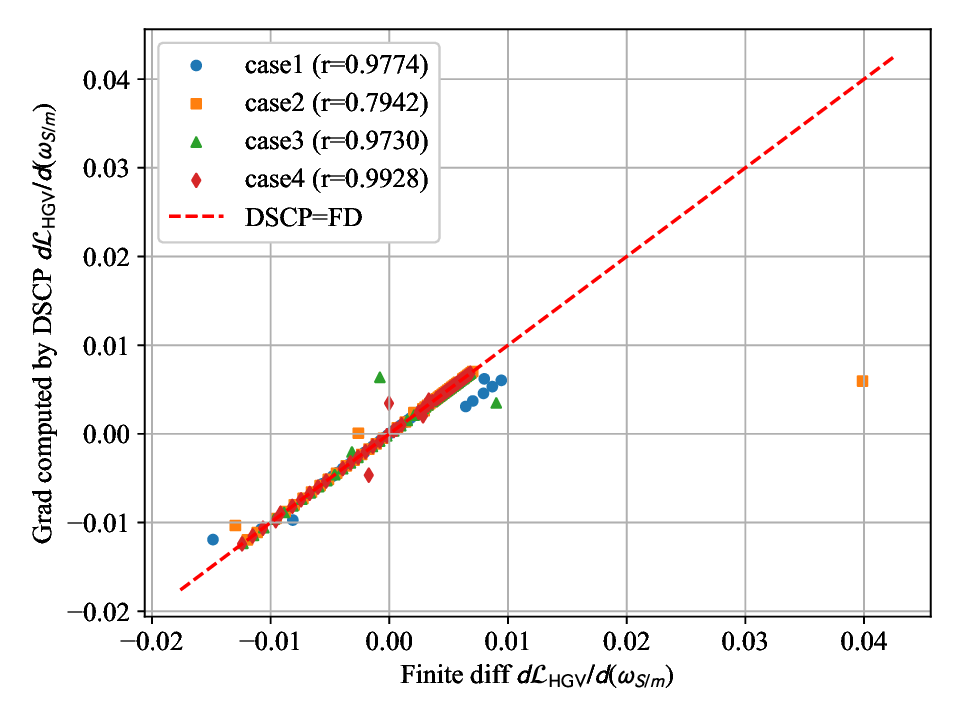}
        \caption{Comparison of gradients with $N_{\mathrm{HGV}}=100$.}
        \label{fig:331}
    \end{subfigure}
    \hfill
    \begin{subfigure}{0.49\textwidth}
        \centering
        \includegraphics[width=\textwidth]{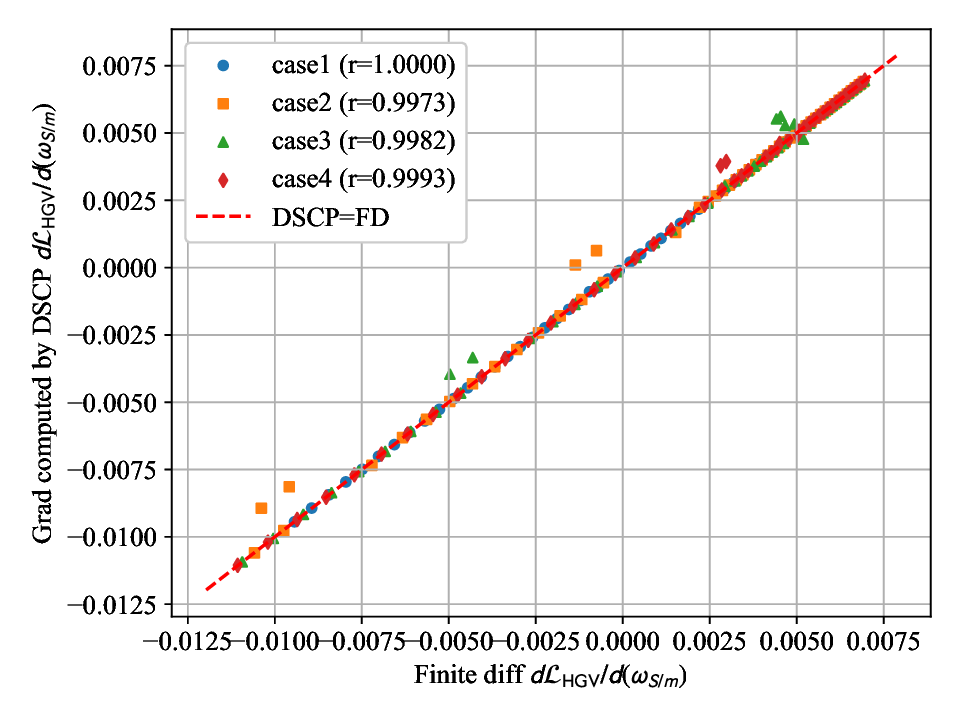}
        \caption{Comparison of gradients with $N_{\mathrm{HGV}}=200$.}
        \label{fig:332}
    \end{subfigure}
    \caption{Comparison of DSCP-computed gradients and finite-difference gradients}
    \label{fig:33}
\end{figure*}

Figures~\ref{fig:31}-~\ref{fig:33} illustrate the overall optimization behavior. Figure~\ref{fig:311} shows the performance index as a function of $ \omega_{S/m}$ for $ N_{\mathrm{HGV}}=100 $. Although the curves corresponding to different trust-region parameters do not coincide exactly, they exhibit highly similar shapes and yield nearly identical optimal values of $\omega_{S/m}$. Notably, several nonsmooth points are observed along the curves, at which the derivative of the cost function is discontinuous. As shown in Fig.~\ref{fig:312}, these nonsmooth locations coincide precisely with changes in the number of SCP iterations required for convergence, indicating that the observed nonsmoothness originates from the numerical behavior of the SCP procedure. Since SCP approximates a continuous-time optimal-control problem through both temporal discretization and iterative convexification, abrupt variations in the number of required iterations may occur. This nonsmooth behavior is progressively mitigated as the convergence tolerance is tightened and the number of discretization nodes is increased. This trend is confirmed by the results in Fig.~\ref{fig:32} for $N_{\mathrm{HGV}}=200$, where the cost function becomes significantly smootherand the iteration numbers exhibit improved stability.

The statistical results in Table~\ref{tab:performance} indicate that the optimal surface-to-mass ratios for both discretization intervals are around $\omega_{S/m}^{\star}=0.92$. This suggests that the nominal airframe design is close to optimal with respect to the chosen performance index, and deviations from optimality may stem from additional geometric or aerodynamic factors that are not explicitly modeled in this study. Moreover, the trust-region coefficients optimized in Sec.~\ref{ch4-2} continue to exhibit favorable numerical performance, with Case 4 achieving the smallest iteration numbers across all experiments. 

Ideally, the gradients computed via DSCP should closely match the finite-difference approximation, and align along the line $ y=x $. Figure~\ref{fig:33} compares the two gradient evaluations. For all four trust-region configurations, the gradient pairs cluster tightly around the diagonal line, with correlation coefficients typically exceeding 0.97 and averaging approximately 0.99. Only a noticable deviation is observed in Case 2 for $ N=100 $. Furthermore, Fig.~\ref{fig:332} shows that the correlation improves substantially when the discretization intervals is increased to $N_{\mathrm{HGV}}=200$, indicating that denser discretization enhances gradient accuracy.

Although a few gradient samples fall in the second or fourth quadrants, indicating occasional sign mismatches, their occurrence is extremely small. Overall, the gradients computed using DSCP exhibit sufficient accuracy and correct directional information to support effective parameter optimization. These results demonstrate that the proposed DSCP framework is capable of providing reliable gradients even for complex SCP formulations and further confirm its applicability to the aerodynamic and geometric design of HGVs.

\section{Conclusions}
\label{ch5}

This paper develops the differentiable sequential convex programming (DSCP) framework for learning and optimizing trajectory planning parameters in aerospace applications. By leveraging differentiable convex optimization, the proposed approach derives sensitivity relations between problem parameters and the optimal solution of SOCP based on first-order optimality conditions. These sensitivities are then combined with the chain rule and gradient backpropagation, enabling end-to-end differentiability from the optimization objective to the underlying parameters. As a result, differentiable mappings are established that support the pretraining of decision variables, the tuning of hyperparameters in trajectory optimization and SCP formulations, and parameter learning within trajectory planning pipelines. The proposed framework is evaluated in two representative scenarios: powered descent guidance and HGV entry. Optimization tasks include the terminal time optimization in powered descent, the trust-region penalty parameters tuning in SCP subproblems, and the surface-to-mass ratio optimization for HGVs. Simulation results demonstrate that the DSCP framework can reliably perform parameter optimization and prelearning of optimization variables. These findings indicate that the proposed approach provides a powerful and versatile tool for vehicle design, mission design, and hyperparameter selection in trajectory planning algorithms across a broad range of aerospace applications.

\section*{Appendix}
\subsection{Details in Powered Descent Guidance}
\label{appendix:PDG}
In this subsection, the specific settings and constructions of subproblems for PDG are provided. For convenience, the subscripts "PDG" are omitted in the following equations. The subproblem for PDG is expressed as
\begin{equation}
    \begin{aligned}
        \min \quad & -m[N]+\omega_{\Gamma}\sum_{n=0}^{N-1}\left\|\Gamma[n]\right\|_2^2\\
        \mathrm{s.t.} \quad & \mathrm{boundary\,constraints}:\\
        & \boldsymbol{x}[0]=\boldsymbol{x}_{0}, \boldsymbol{r}[N]=\boldsymbol{r}_{\mathrm{f}},\boldsymbol{v}[N]=\boldsymbol{v}_{\mathrm{f}}\\
        & \mathrm{discretized\,dynamics\,constraints}\,(n\in [0,N-1]):\\
        & \left(\frac{T_{\text{ref}}}{2}\left.\frac{\partial \boldsymbol{f}}{\partial \boldsymbol{x}}\right|_{\text{ref}} + \boldsymbol{I}\right)\delta \boldsymbol{x}[n] + \left(\frac{T_{\text{ref}}}{2}\left.\frac{\partial \boldsymbol{f}}{\partial \boldsymbol{x}}\right|_{\text{ref}} - \boldsymbol{I}\right)\delta \boldsymbol{x}[n+1]+ \frac{T_{\text{ref}}}{2}\left.\frac{\partial \boldsymbol{f}}{\partial \boldsymbol{u}}\right|_{\text{ref}}\delta  \boldsymbol{u}[n] + \frac{T_{\text{ref}}}{2}\left.\frac{\partial \boldsymbol{f}}{\partial \boldsymbol{u}}\right|_{\text{ref}} \delta \boldsymbol{u}[n+1]\\
        &=(\boldsymbol{x}_{\text{ref}}[n+1]- \boldsymbol{x}_{\text{ref}}[n])-T_{\text{ref}} \boldsymbol{f}_{\text{ref}}[n]+ \boldsymbol{\Gamma}[n]\\
        & \mathrm{thrust\,constraints}\,(n\in [0,N]):\\
        & u_{\min}\leq \left\|\boldsymbol{u}_{\mathrm{ref}}[n]\right\|_{2}+\frac{\boldsymbol{u}_{\mathrm{ref}}^{\top}[n]}{\left\|\boldsymbol{u}_{\mathrm{ref}}[n]\right\|_{2}} \left(\boldsymbol{u}[n]-\boldsymbol{u}_{\mathrm{ref}}[n]\right) \leq u_{\max}\\
        &\sqrt{u_y^2[n]+u_z^2[n]}\leq \tan(\eta_{\max})u_x[n]\\
        & \mathrm{path\,constraints}\,(n\in [0,N]):\\
        & \sqrt{r_y^2[n]+r_z^2[n]}\leq \tan(\beta_{\max}) r_x[n]\\
    \end{aligned}
\end{equation}
where $\boldsymbol{\Gamma}$ is the virtual control vector and $\omega_{\Gamma}=10000$ is the weighting coefficient for the virtual control. The dynamics Jacobians are given by
\begin{equation}
    \frac{\partial \boldsymbol{f}}{\partial \boldsymbol{x}}=\begin{bmatrix}
    \boldsymbol{0}_{3\times 3} & \boldsymbol{I}_{3\times 3} & \boldsymbol{0}_{3\times 1} \\
    \boldsymbol{0}_{3\times 3} & \boldsymbol{0}_{3\times 3} & -\frac{\boldsymbol{u}}{m} \\
    \boldsymbol{0}_{1\times 3} & \boldsymbol{0}_{1\times 3} & 0
    \end{bmatrix}, \quad
    \frac{\partial \boldsymbol{f}}{\partial \boldsymbol{u}}=\begin{bmatrix}
    \boldsymbol{0}_{3\times 3} \\ \frac{\boldsymbol{I}_{3\times 3}}{m} \\ -\frac{\boldsymbol{u}^{\top}}{I_{\mathrm{sp}}g_0\left\|\boldsymbol{u}\right\|_2}
\end{bmatrix}
\end{equation}

The initial guess of the SCP procedure for PDG is provided by a simplified problem given by 
\begin{equation}
    \begin{aligned}
        \min \quad & \sum_{n=0}^{N}\left\|\boldsymbol{u}[n]\right\|_2^2\\
        \mathrm{s.t.} \quad
        & \boldsymbol{x}[0]=\boldsymbol{x}_{0}, \boldsymbol{r}[N]=\boldsymbol{r}_{\mathrm{f}},\boldsymbol{v}[N]=\boldsymbol{v}_{\mathrm{f}}\\
        & \boldsymbol{r}[n+1]-\boldsymbol{r}[n]=\frac{T}{2}\left(\boldsymbol{v}[n+1]+\boldsymbol{v}[n]\right)\\
        & \boldsymbol{v}[n+1]-\boldsymbol{v}[n]=\frac{T}{2}\left(\frac{\boldsymbol{u}[n+1]}{m[n+1]}+\frac{\boldsymbol{u}[n]}{m[n]}+2\boldsymbol{g}_{E}\right),m[n]=m_{0}\\
    \end{aligned}
\end{equation}
In the above problem, the variation of mass and the thrust constraints are disregarded, so the resulting trajectory does not satisfy the constraints of PDG. In the early stages of the iteration, the virtual control can effectively compensate for the deficiencies of the initial trajectory.

\subsection{Details in Entry Trajectory Optimization}
\label{appendix:HGV}
In this subsection, we detail the formulation and setup of the subproblems associated with entry trajectory planning for HGVs. For simplicity, the subscript "HGV" is omitted in the subsequent equations. The subproblem can be written as
\begin{equation}
    \begin{aligned}
        \min \quad & \mathrm{Eq}.\eqref{equation:j2hgv}\\
        \mathrm{s.t.} \quad & \mathrm{boundary\,constraints}:\\
        & \boldsymbol{x}[0]=\boldsymbol{x}_{0}, \boldsymbol{x}[N]=\boldsymbol{x}_{\mathrm{f}}(\mathrm{except\,for},V[N]),V[n]\geq V_{\mathrm{f,\min}}\\
        & \mathrm{discretized\,dynamics\,constraints}\,(n\in [0,N-1]):\\
        & \left(\frac{T_{\text{ref}}}{2}\left.\frac{\partial \boldsymbol{f}}{\partial \boldsymbol{x}}\right|_{\text{ref}} + \boldsymbol{I}\right)\delta \boldsymbol{x}[n] + \left(\frac{T_{\text{ref}}}{2}\left.\frac{\partial \boldsymbol{f}}{\partial \boldsymbol{x}}\right|_{\text{ref}} - \boldsymbol{I}\right)\delta \boldsymbol{x}[n+1]+ \frac{T_{\text{ref}}}{2}\left.\frac{\partial \boldsymbol{f}}{\partial \boldsymbol{u}}\right|_{\text{ref}}\delta  \boldsymbol{u}[n] + \frac{T_{\text{ref}}}{2}\left.\frac{\partial \boldsymbol{f}}{\partial \boldsymbol{u}}\right|_{\text{ref}} \delta \boldsymbol{u}[n+1]\\
        &=(\boldsymbol{x}_{\text{ref}}[n+1]- \boldsymbol{x}_{\text{ref}}[n])-T_{\text{ref}} \boldsymbol{f}_{\text{ref}}[n]\\
        & \mathrm{path\,constraints}\,(n\in [0,N]):\\
        & \bar{p}[n]+\left.\frac{\partial \bar{p}}{\partial r}\right|_{\text{ref}}\delta r[n]+\left.\frac{\partial \bar{p}}{\partial V}\right|_{\text{ref}}\delta V[n]\leq \bar{p}_{\max}\\
        & \bar{q}[n]+\left.\frac{\partial \bar{q}}{\partial r}\right|_{\text{ref}}\delta r[n]+\left.\frac{\partial \bar{q}}{\partial V}\right|_{\text{ref}}\delta V[n]\leq \bar{q}_{\max}\\
        & \bar{n}[n]+\left.\frac{\partial \bar{n}}{\partial r}\right|_{\text{ref}}\delta r[n]+\left.\frac{\partial \bar{n}}{\partial V}\right|_{\text{ref}}\delta V[n]+\left.\frac{\partial \bar{n}}{\partial \alpha}\right|_{\text{ref}}\delta \alpha[n]\leq \bar{n}_{\max}\\
        & \mathrm{control\,angle\,constraints}\,(n\in [0,N]):\\
        & \alpha_{\min}\leq \alpha[n]\leq \alpha_{\max},\, \sigma_{\min}\leq \sigma[n]\leq \sigma_{\max}\\
        & \dot\alpha_{\min}\leq \dot\alpha[n]\leq \dot\alpha_{\max},\, \dot\sigma_{\min}\leq \dot\sigma[n]\leq \dot\sigma_{\max}\\
    \end{aligned}
\end{equation}
The dynamics Jacobians and the state-dependent differentials of the path constraints are omitted here for brevity and can be found in \cite{wang2017constrained}. The initial trajectory for all entry simulations is calculated using a fixed control profile, which is given by
\begin{equation}
    \begin{aligned}
    &u_1(t)= \frac{\frac{\pi}{9} - \alpha_0}{t_{\mathrm{f0}}}\\
    &u_2(t)=-\,\frac{2\pi}{3}\,\frac{\sigma_{\max}}{t_{\mathrm{f0}}}
    \cos\!\left(\frac{2\pi t}{t_{\mathrm{f0}}}\right)
  - \frac{1}{6}\,\frac{\sigma_{\max}}{t_{\mathrm{f0}}}.
    \end{aligned}
\end{equation}
where $t_{\mathrm{f0}}=2500s$ is the final time for the initial reference trajectory, $\alpha_0$ is the initial angle of attack, and $\sigma_{\max}$ is the bank angles, respectively. These settings guarantee that the initial reference controls satisfy the angle and angular-rate constraints in Eq.\eqref{equation:angleconstraint}.

\section*{Acknowledgments}
This work was supported by the National Natural Science Foundation of China (Grant No. 12525204) and the Natural Science Foundation of Beijing Municipality (Grant No. L251014).

\bibliography{sample}

\end{document}